%% file: main.tex
\documentclass[11pt]{article}
\usepackage[margin=1.1in]{geometry}
\usepackage{amsmath,amssymb,amsthm,mathtools}
\usepackage{enumitem}
\usepackage{microtype}
\usepackage{xcolor}
\definecolor{linkblue}{rgb}{0.05,0.2,0.55}
\usepackage[colorlinks=true,linkcolor=linkblue,citecolor=linkblue,urlcolor=linkblue]{hyperref}

\theoremstyle{plain}
\newtheorem{theorem}{Theorem}[section]
\newtheorem{lemma}[theorem]{Lemma}
\newtheorem{proposition}[theorem]{Proposition}
\newtheorem{corollary}[theorem]{Corollary}
\theoremstyle{remark}
\newtheorem{remark}[theorem]{Remark}
\numberwithin{equation}{section}

\newcommand{\T}{\mathbb{T}}

\title{The sharp uniqueness threshold for the\\ antiferromagnetic Potts model on regular trees}
\date{\today}

\author{%
Wei Wu\thanks{(W. Wu) Department of Mathematics, NYU Shanghai, and NYU-ECNU Institute of Mathematical Sciences at NYU Shanghai, China.
\texttt{wei.wu@nyu.edu}}
\and
Yicheng Xu\thanks{(Y. Xu) Department of Mathematics, NYU Shanghai, China.
\texttt{yx4167@nyu.edu}}
\and
Kuan Yang\thanks{(K. Yang) John Hopcroft Center for Computer Science, Shanghai Jiao Tong University, China.
\texttt{kuan.yang@sjtu.edu.cn}}%
}

\begin{document}
\nocite{*}
\maketitle

\begin{abstract}
We study Gibbs uniqueness for the $q$-state antiferromagnetic Potts
model on the infinite rooted $d$-ary tree and the infinite
$(d+1)$-regular tree, where $d\geqslant1$ and $q\geqslant2$.
An edge has weight $p\in[0,1]$ if its endpoints have the same color
and weight $1$ otherwise.
We prove that there is a unique Gibbs state if and only if
$p\geqslant p_c:=\max\{0,1-q/(d+1)\}$,
except when $(q,p)=(d+1,0)$, thus completes the uniqueness thresholds for all $d\geqslant1$ and $q\geqslant2$.
Our proof use a convex potential
for the tree recursion and show that it strictly decreases
away from the uniform distribution.
Together with the result of Gu, Wu and Yang, our theorem
gives the exact critical power law whenever $p_c>0$.
\end{abstract}

\tableofcontents

\input{sections/01-introduction.tex}
\input{sections/02-recursion.tex}
\input{sections/03-logarithmic-comparison.tex}

\input{sections/04-scalar.tex}
\input{sections/05-uniform-convergence.tex}

\input{sections/06-gibbs-uniqueness.tex}

\bibliographystyle{manuscript}
\bibliography{refs}

\end{document}

%% file: sections/01-introduction.tex
\section{Introduction}\label{sec:intro}

In this paper, we study the uniqueness of Gibbs measures for the antiferromagnetic (AF) $q$-state Potts model on regular trees. We are interested in the influence of a boundary condition at distance $n$ on the distribution of the spin at root. Our main result determines the exact range of the edge parameter for which this influence vanishes as $n$ tends to infinity. In particular, we establish the sharp uniqueness threshold for the AF Potts model.

Let $\mathbb T^d$ be the infinite rooted $d$-ary tree, where $d\geqslant1$. We denote its root by $\ast$ and the subtree of height $n$ by
$\mathbb T^d_n$. Its boundary $\partial\mathbb T^d_n$ consists of the vertices at distance $n$ from the root, and we write
$\operatorname{int}(\mathbb T^d_n)$ for its interior. Following the dictionary notation for the tree in~\cite{GWY}, we identify a vertex $u$ at distance $k$ with a word of length $|u|=k$ over $\{1,\ldots,d\}$. The root is an empty word, and the children of $u$ are $u1,\ldots,ud$.
Then we can write \(\T^d_n := \cup_{k=0}^n \{1,\cdots,d\}^k \) and $\T^d := \cup_{k=0}^\infty \{1,\cdots,d\}^k$.
We use the notation $u\sim v$ when two vertices are adjacent.

Let $q\geqslant2$ be the number of colors, and write $[q]=\{1,\ldots,q\}$. For a spin configuration $\sigma:\mathbb T^d_n\to[q]$, let $M(\sigma)$ denote the number of monochromatic edges, namely, edges whose endpoints have the same color. Given a boundary condition $\xi:\partial\mathbb T^d_n\to[q]$ and an edge parameter $p\in[0,1]$, the weight of $\sigma$ is given by
\begin{equation}\label{eq:W}
W^\xi_{n,p,q}(\sigma)
=p^{M(\sigma)}\mathbf{1}_{\{\sigma|_{\partial\mathbb T^d_n}=\xi\}}.
\end{equation}
Normalizing these weights gives the Potts measure on $\mathbb T^d_n$
with boundary condition $\xi$,
\begin{equation}\label{eq:mu}
\mu^\xi_{n,p,q}[\sigma]=\frac{W^\xi_{n,p,q}(\sigma)}{Z^\xi_{n,p,q}},
\qquad
Z^\xi_{n,p,q}=\sum_{\sigma\in[q]^{\mathbb T^d_n}}W^\xi_{n,p,q}(\sigma).
\end{equation}
For $0<p<1$, the interaction favors different colors at adjacent vertices, since every
monochromatic edge contributes a factor $p$ to the weight. At zero
temperature ($p=0$), this becomes the uniform distribution over proper $q$-colorings that extend $\xi$ (with the convention $0^0=1$). We say that a boundary condition $\xi$ is admissible if $Z^\xi_{n,p,q}>0$. Note that all boundary conditions
are admissible if $p>0$. Throughout the paper, extremal values are taken over admissible boundary conditions, and $d,q,p$ are fixed as $n$ tends to infinity.

On the infinite tree we use the Dobrushin--Lanford--Ruelle definition of Gibbs measures. For each finite set of vertices, conditional on the spins outside the set, the spins inside have the Potts distribution with these exterior values. At $p=0$ the measure is supported on all proper colorings.

A fundamental question in statistical physics is whether the interaction at
a given temperature determines a single equilibrium state, which is also known as the uniqueness/nonuniqueness of the Gibbs measure. In the case
of uniqueness, the marginal distribution of the spins in any fixed finite region
has the same limit as the distance to the boundary increases, no matter what colors are fixed
there. Coexistence of Gibbs measures means that the choice of distant
spins can still be detected locally.
A change from uniqueness to nonuniqueness as the temperature varies marks a phase transition.

For the Potts model on regular trees, uniqueness can be characterized
by the convergence of the marginal distribution at the root
\cite[Lemma~1.3]{dBBR}. This marginal must converge to the uniform
distribution as $n$ tends to infinity, uniformly over the spins fixed
on $\partial\mathbb T^d_n$. By symmetry, it suffices to consider color~$1$,
and the criterion becomes
\begin{equation}\label{eq:Delta}
\lim_{n\to\infty}\max_\xi
\left|\mu^\xi_{n,p,q}[\sigma_\ast=1]-\frac1q\right|=0.
\end{equation}

It is a folklore conjecture that there exists a critical parameter
\begin{equation}\label{eq:pc}
p_c=p_c(d,q):=\max\left\{0,\,1-\frac{q}{d+1}\right\}
\end{equation}
such that the Gibbs measure of the Potts model on the infinite
$d$-ary tree (or the infinite $(d+1)$-regular tree) is unique if and only if
$p\in[p_c,1]$ and $(q,p)\ne(d+1,0)$.
To see how the critical parameter arises, consider the case where each
of the $d$ subtrees has root marginal proportional to
$(r,1,\ldots,1)$. The ratio of the probabilities of colors $1$ and $2$
is $r$ at each child. An iteration of the tree recursion gives the
corresponding ratio in the parent as $f(r)$, where
\[
f(r)=\left(\frac{pr+q-1}{r+p+q-2}\right)^d,
\qquad f(1)=1,\qquad f'(1)=-\frac{d(1-p)}{p+q-1}.
\]
The fixed point $r=1$ represents the uniform distribution. It is locally
attracting when $|f'(1)|<1$, that is, when $p>1-q/(d+1)$. This suggests
the critical parameter $p_c$.
For an arbitrary boundary condition, the marginals of the $d$ subtrees do not need to coincide, so their influence on the root cannot in
general be described by a single ratio $r$, which yields a much more complicated recursion.  There is a further difficulty
at $p=p_c>0$. Even when the child marginals have the form above,
$f'(1)=-1$, and the linear approximation does not decay. The conjecture
is that \eqref{eq:Delta} holds for all $p\in[p_c,1]$, with the exception of
$p=0$ when $q=d+1$ (see \cite{GWY}).

The conjecture had been established in several cases.
At $p=0$, which is the proper $q$-coloring model, Jonasson~\cite{Jonasson}
proved the uniqueness for $q\geqslant d+2$. When $q\leqslant d+1$, frozen
colorings give different Gibbs measures~\cite{BW}. At positive temperatures,
it has been known that nonuniqueness occurs as long as $p < p_c$~\cite{PLM83}.
The remaining question was whether uniqueness holds down to $p_c$ for arbitrary $q$ and $d$, including at the critical point $p_c$ itself whenever it is positive.
Galanis, Goldberg and Yang~\cite{GGY}
confirmed the conjecture for $q=3$, and $q=d=4$, and de Boer, Buys and
Regts~\cite{dBBR} proved it for $q=4$ and $d\geqslant4$. For larger
$q$, Bencs, de Boer, Buys and Regts~\cite{BdBBR} showed that uniqueness
holds for $p\geqslant p_c$ for each fixed $q\geqslant 5$ and all
sufficiently large $d$.
Further progress towards the conjectured threshold was made by Chen, Liu, Mani, and Moitra~\cite{CLMM}, who proved uniqueness for $p\geq 1-\frac{q-1}{d+6}$, and subsequently by Bencs, Berrekkal and Regts~\cite{BBR}, who established uniqueness down to $p_c+O(q/d^2)$ for $q\geq 3$ and $d\geq q+2$.
However, none of them completely solved the conjecture. In this work, we close the remaining gap and establish the conjectured uniqueness threshold.

The uniqueness threshold and the influence of the boundary also play an important role in the complexity of
counting and sampling on finite graphs. For the AF Ising model without external field ($q=2$), fix an integer $\Delta\geqslant3$.
The tree uniqueness threshold separates two regimes for approximating the partition function on bounded degree graphs. For $p>p_c=1-2/\Delta$, Sinclair, Srivastava
and Thurley~\cite{SST} gave a deterministic algorithm on graphs of maximum degree $\Delta$. It approximates the partition function to relative error $\varepsilon$ in time polynomial in the number of vertices and $1/\varepsilon$. For $0<p<p_c$, Sly and Sun~\cite{SS} showed that there is no randomized approximation scheme on $\Delta$-regular graphs unless $\mathrm{NP}=\mathrm{RP}$. At $p=p_c$, Chen, Chen, Yin and Zhang~\cite{CCYZ} proved that single-site Glauber dynamics mixes in polynomial time on graphs of maximum degree $\Delta$, giving a sampling algorithm at the critical value. For $q\geqslant3$, it is also a folklore conjecture that there exists a polynomial-time approximation schemes (and a polynomial-time sampler) for the Potts model on graphs of maximum degree $\Delta$, when parameters are in the uniqueness regime. This conjecture is partially solved on random instances.
Blanca et al.~\cite[Theorem~2.7]{BGGSVY} gives a sampling algorithm for random regular graphs, depending on the uniqueness of Gibbs measure on the tree. More precisely, fix $q,\Delta\geqslant3$ and $0<p<1$ with
$p\neq1-q/\Delta$. If uniqueness holds on the $(\Delta-1)$-ary tree,
then with probability tending to one over the choice of the random
$n$-vertex $\Delta$-regular graph, the algorithm produces in polynomial
time a sample whose law is within $O(n^{-\delta})$ of the Potts distribution
in total variation, for some $\delta=\delta(q,\Delta,p)>0$. Theorem~\ref{thm:main}
below verifies their uniqueness hypothesis for every $p_c<p<1$. Whether
efficient sampling is possible throughout $p_c<p<1$ on all graphs of
maximum degree $\Delta$ remains open for general $q$.
Below the tree threshold, Galanis, \v{S}tefankovi\v{c} and
Vigoda~\cite{GSV} proved that for even $q\geqslant4$ and
$0\leqslant p<1-q/\Delta$, there is no fully polynomial randomized
approximation scheme for the partition function, even on triangle-free
$\Delta$-regular graphs, unless $\mathrm{NP}=\mathrm{RP}$.

\subsection{Main result}

Our main theorem establishes the conjectured uniqueness threshold on both infinite $d$-ary tree and $(d+1)$-regular tree.

\begin{theorem}[Sharp uniqueness threshold]\label{thm:main}
Let $q,d\in\mathbb{N}_+$ with $q\geqslant 2$, and let $p\in[0,1]$. The $q$-state AF Potts model on $\mathbb{T}^d$ has a unique Gibbs measure if and only if
\begin{equation}\label{eq:range}
\begin{cases}
p\in[p_c,1], & q\neq d+1,\\
p\in(p_c,1], & q=d+1.
\end{cases}
\end{equation}
The same criterion holds on the infinite $(d+1)$-regular tree.

In this regime, the root marginal satisfies
\begin{equation}\label{eq:Delta}
\lim_{n\to\infty}
\max_{\xi}
\left|
\mu^\xi_{n,p,q}[\sigma_\ast=1]-\frac{1}{q}
\right|
=0.
\end{equation}
More generally, let $\nu$ be the unique Gibbs measure on $\mathbb{T}^d$. For every fixed $m\geqslant 0$ and every event $A$ depending only on the spins in $\mathbb{T}^d_m$,
\begin{equation}\label{eq:unifconv}
\lim_{n\to\infty}
\sup_{\xi}
\left|
\mu^\xi_{n,p,q}[A]-\nu[A]
\right|
=0.
\end{equation}
Outside \eqref{eq:range}, each tree admits at least two Gibbs measures, and the convergence in \eqref{eq:Delta} fails.
\end{theorem}

The unique measure in Theorem~\ref{thm:main} is the free
Gibbs measure. On either tree, it is obtained by assigning
a uniform spin to the root and then sampling the spins
recursively, independently across children conditional on
their parent's spin, with transition probabilities
\begin{equation}\label{eq:Q}
K_p(j\mid i):=
\frac{1+(p-1)\delta_{i,j}}{p+q-1},
\qquad i,j\in[q],
\end{equation}
where $\delta_{i,j}$ is the Kronecker delta.

The exception at $q=d+1$ occurs at zero temperature. A proper coloring can be frozen by its boundary, so that the color of the root is determined however far away the boundary is. For the same $q$ and $d$, Theorem~\ref{thm:main} gives uniqueness for every $p>0$.

For $2\leqslant q\leqslant d$, the critical parameter is positive,
and one can also ask how fast the influence of the boundary vanishes
at $p_c$. Gu, Wu and Yang~\cite[Theorem~1.2]{GWY} determined the exact decay rate under the assumption of uniqueness. Their estimate applies once the root marginal is sufficiently close to the uniform distribution.
Theorem~\ref{thm:main} ensures that every boundary condition gives a marginal in this neighborhood for all sufficiently large $n$. Combining the two results yields the following power law.

\begin{corollary}[Power law decay]\label{cor:powerlaw}
For $2\leqslant q\leqslant d$ and $p_c=1-q/(d+1)$,
\begin{equation}\label{eq:powerlaw}
\lim_{n\to\infty}\sqrt n\,
\max_{\xi:\partial\mathbb T^d_n\to[q]}
\left|\mu^\xi_{n,p_c,q}[\sigma_\ast=1]-\frac1q\right|
=\sqrt{\frac{6d^2}{d^2-1}}\,\frac{q-1}{q^2}.
\end{equation}
\end{corollary}


\subsection{Outline of the proof}\label{sec:strategy}

The main difficulty is to control the influence of arbitrary boundary conditions. The root marginals of different child subtrees need not coincide, so the recursion must allow different distributions on these subtrees. We handle this by constructing an explicit convex potential function for every $d\geqslant2$ and $2\leqslant q\leqslant d+1$, that decreases globally under a one-step recursion at every nonzero centered log-message.

We first write a tree recursion in logarithmic coordinates to compute marginals. 
Let $x_i^{(u)}$ denote the probability of color $i$ at $u$ in
this subtree, with the edge to the parent removed. If the parent has
color $i$, adding the edge gives weight $p$, otherwise it gives weight $1$. Averaging over the spin at $u$ therefore gives
the factor $px_i^{(u)}+1-x_i^{(u)}=1-(1-p)x_i^{(u)}$. Multiplying
the contributions of the $d$ subtrees, we obtain for positive child
marginals the following identity for the ratio at the parent
\[
\frac{x'_i}{x'_j}=\prod_{u=1}^d
\frac{1-(1-p)x^{(u)}_i}{1-(1-p)x^{(u)}_j}.
\]
Observe that the $d$ vectors on the right depend on the restrictions
of $\xi$ to the respective subtrees. To compare them we use centered
logarithmic coordinates. For a strictly positive marginal $x$, let
\[
y_i:= \log x_i-\frac1q\sum_{j=1}^q\log x_j,
\]
Then $y_i-y_j=\log(x_i/x_j)$ for each pair of colors. 
Note that $\sum_i x_i = 1$. Thus we can uniquely recover $x$ from $y$ by
\[
    x_i = \frac{e^{y_i}}{\sum_j e^{y_j}}\,.
\]
In particular, $y = 0$ corresponds to the uniform distribution.
Moreover, permuting the colors simply permutes the coordinates of $y$. 

In our tree recursion, the contributions of the child subtrees of each node multiply. Taking logarithms turns this product into a sum, and subtracting the mean removes the common normalizing factor.
Section~\ref{sec:messages} shows that the recursion takes the form
\begin{equation}\label{eq:avg-intro}
y'=\frac1d\sum_{u=1}^d T_{d,p}\bigl(y^{(u)}\bigr),
\end{equation}
where $y^{(u)}$ is the centered logarithmic vector for the root marginal of the subtree at child $u$, $y'$ is the corresponding vector at the parent node, and the map $T_{d,p}$ describes the case in which all child marginals coincide. In reference color log-ratio coordinates, this averaging identity appears in the work of Bencs, de Boer, Buys and Regts\cite[Lemma~2.2]{BdBBR}.

The coordinates of $y$ sum to zero. We therefore work on the space
\begin{equation}\label{eq:Eq}
E_q:=\left\{y\in\mathbb{R}^q:\sum_{i=1}^q y_i=0\right\},
\end{equation}
and consider a potential function
\begin{equation}\label{eq:Phic}
\Phi_c(y):=\sum_{i=1}^q e^{cy_i}-q,
\qquad \text{for some } c>0.
\end{equation}
The function $\Phi_c$ is convex and nonnegative on $E_q$, and vanishes only at $0$. We will choose $c$ so that this potential function decreases under the tree recursion for computing marginals.
The child marginals in \eqref{eq:avg-intro} may be different. By the convexity of $\Phi_c$, Jensen's inequality gives
\[
\Phi_c(y')
\leqslant
\frac1d\sum_{u=1}^d
\Phi_c\bigl(T_{d,p}(y^{(u)})\bigr).
\]
Thus, if $\Phi_c(T_{d,p}(y))<\Phi_c(y)$ for every nonzero $y\in E_q$,
the potential at the parent is smaller than the largest child
potential, unless all child vectors are zero. This is how an
inequality for $T_{d,p}$ applies to arbitrary child marginals.

The exponential form of $\Phi_c$ also yields
\[
e^{cy_i}
=
\frac{x_i^c}{\left(\prod_{j=1}^q x_j\right)^{c/q}}.
\]
We use this identity, and the corresponding identity for
$T_{d,p}(y)$, to compare the two potentials in Section~\ref{sec:entropy}.
Moreover, a bound on $\Phi_c$ bounds each coordinate of $y$ from above.
Since the coordinates sum to zero, it also bounds them from below.
This gives the compactness property stated in Lemma~\ref{lem:potential}.

Now consider $d\geqslant2$ and $2\leqslant q\leqslant d+1$.
At $p=p_c$, the derivative of $T_{d,p_c}$ at $0$ is $-\mathrm{Id}$,
so the linear approximation gives no strict contraction.
To see how to choose the exponent $c$, let $\gamma>1$ and define
\[
F_{q,d,\gamma}(x)
:=
(qx)^\gamma
-
\left(\frac{d+1-qx}{d}\right)^{\gamma d}
-\gamma\log(qx)
+\gamma d\log\left(\frac{d+1-qx}{d}\right),
\qquad
0<x<\frac{d+1}{q}.
\]
In Section~\ref{sec:entropy}, we will see that it suffices  to prove 
\begin{equation}\label{eq:scalar-intro}
    \sum_{i=1}^q F_{q,d,\gamma}(x_i)>0,
\end{equation}
for every positive nonuniform probability vector $x$.

We first examine this condition close to the uniform marginal. For
any $\gamma >1$, the value of $F_{q,d,\gamma}$ and its first two derivatives at $1/q$ vanish. Then next terms of its Taylor expansion gives
\[
F_{q,d,\gamma}\left(\frac1q+h\right)
=
\frac{q^3\gamma^2}{3}
\left(\gamma-\frac{3(d+1)}{2d}\right)h^3
+O(h^4).
\]
A nonzero cubic term would give opposite signs on the two sides of
$1/q$ for sufficiently small $|h|$. To satisfy this local requirement,
we must cancel that term. We therefore take
\[
\gamma=c_d:=\frac{3(d+1)}{2d}.
\]
At this value,
\[
F_{q,d,\gamma}^{(4)}(1/q)
=
2q^4\gamma^2\left(1-\frac1{d^2}\right)>0,
\qquad
F_{q,d,\gamma}(1/q+h)
=
\frac{q^4\gamma^2}{12}\left(1-\frac1{d^2}\right)h^4+O(h^5).
\]
Thus the selected scalar function is positive on both sides of $1/q$
in a sufficiently small neighborhood.
Let $F_{q,d}(x) := F_{q,d,c_d}(x)$ be the function obtained by fixing $\gamma$.

To prove positivity on the whole interval $(0,d/q]$,
we study the shape of $F_{q,d}$ and prove that $F_{q,d}(d/q)>0$
in Section~4. 
Then these estimates imply
\[
F_{q,d}(x)>0
\qquad
\text{for }0<x\leqslant d/q,\quad x\ne1/q.
\]
Since $F_{q,d}(1/q)=0$, the sum is strictly positive whenever all
coordinates satisfy $x_i\leqslant d/q$ and the vector is nonuniform.
This includes every case with $q\leqslant d$.

The remaining case is $q=d+1$ with a coordinate, say $x_1$, larger
than $d/q$. The other $d$ coordinates have sum $1-x_1<\frac1q$,
so each lies in $(0,1/q)$. We prove that $F_{q,d}$ is convex on this
interval. Jensen's inequality then gives
\[
\begin{aligned}
\sum_{i=1}^{d+1}F_{q,d}(x_i)
&\geqslant
F_{q,d}(x_1)
+dF_{q,d}\left(\frac{1-x_1}{d}\right)\\
&=
F_{q,d}(x_1)
+dF_{q,d}\left(\frac{d+1-qx_1}{dq}\right)
=:S_{q,d}(x_1).
\end{aligned}
\]
A separate estimate establishes
\[
S_{q,d}(x)>0
\qquad
\text{for }d/q<x<(d+1)/q.
\]
Thus, even if the large coordinate contributes a negative term,
the remaining coordinates contribute enough to make the sum positive.
This completes the proof of \eqref{eq:scalar-intro} for our choice of $\gamma = c_d$, and hence of the strict potential decrease $\Phi_{c_d}(T_{d,p}(y))<\Phi_{c_d}(y)$ for every nonzero $y\in E_q$.

Returning to the tree recursion, convexity and
\eqref{eq:avg-intro} give
\[
\begin{aligned}
\Phi_{c_d}(y')
&\leqslant
\frac1d\sum_{u=1}^d
\Phi_{c_d}\bigl(T_{d,p_c}(y^{(u)})\bigr) <
\max_{1\leqslant u\leqslant d}
\Phi_{c_d}\bigl(y^{(u)}\bigr),
\end{aligned}
\]
provided at least one child vector is nonzero.

For convergence uniformly over boundary conditions, we need to
bound the number of iterations independently of the boundary.
Suppose all initial potentials are at most some finite $a_0$.
For $0<\varepsilon<a_0$, the set
\[
\left\{y\in E_q:
\varepsilon\leqslant\Phi_{c_d}(y)\leqslant a_0\right\}
\]
is compact. The continuous function
\[
\Phi_{c_d}(y)-\Phi_{c_d}\bigl(T_{d,p_c}(y)\bigr)
\]
therefore has a strictly positive minimum on this set.
Lemma~\ref{lem:compact} uses this fact and the convexity to obtain a number of layers after which every potential is at most $\varepsilon$, where the number depends only on $a_0$ and $\varepsilon$ for fixed $d,q$.
Consequently, a common finite initial bound gives convergence to
the uniform distribution uniformly over boundary conditions.

Section~\ref{sec:uniform} obtains the initial bound after one layer when $p>0$,
and after two layers when $p=0$ and $q\geqslant d+2$.
It also extends the descent to larger $p$ by monotonicity and
handles $q\geqslant d+2$ and $d=1$.
At the exceptional point $(q,p)=(d+1,0)$, frozen boundary conditions
can force zero coordinates, so the logarithmic coordinates need
not be defined.
Finally, Section~\ref{sec:gibbs} deduces Gibbs uniqueness and convergence on fixed finite subtrees, and applies the known nonuniqueness results
outside \eqref{eq:range}.

\paragraph*{Acknowledgments.}
We benefited from the assistance of AI tools in our study. We asked ChatGPT-6 Astra whether it could find a Lyapunov function
for the tree recursion, and used it to write the proofs in Section~\ref{sec:entropy} and Section~\ref{sec:scalar} under our direction. The authors subsequently revised and rewrote the draft. We thank Chenlin Gu for bringing the work \cite{song2026ai} to our attention. The work of W.W. is supported in part by a MOST grant and by the NYUSH-Boost Fund. The work of K.Y. is supported by NSFC grant 62572298.





%% file: sections/02-recursion.tex
\section{The tree recursion and a compactness principle}
\label{sec:messages}

\subsection{Root marginals and the one-step iteration}

Fix integers $d\geqslant1$, $q\geqslant2$ and $n\geqslant1$,
a parameter $p\in[0,1]$, and an admissible boundary condition
$\xi$ on $\partial\mathbb{T}^d_n$.

For $v\in\mathbb{T}^d_n$, the subtree consisting of $v$ and its
descendants has height $h=n-|v|$.
We identify it with $\mathbb{T}^d_h$ by writing each descendant
as $vw$, where $vw$ denotes concatenation.
The boundary condition inherited from $\xi$ is
\[
\xi^v(w):=\xi(vw),
\qquad
w\in\partial\mathbb{T}^d_h.
\]
Define
\begin{equation}\label{eq:Zv}
Z_v:=Z^{\xi^v}_{h,p,q},
\qquad
x_i^{(v)}:=\mu^{\xi^v}_{h,p,q}[\sigma_\ast=i],
\qquad i\in[q].
\end{equation}
Thus $x^{(v)}=(x_i^{(v)})_{i\in[q]}$ is the root marginal of this
subtree, with the edge to the parent of $v$ removed when $v\ne\ast$.
Admissibility of $\xi$ ensures that $Z_v>0$, so these probabilities
are well defined.

When $p>0$, every configuration extending the boundary condition has positive weight, so $Z_v>0$. When $p=0$, admissibility of $\xi$ provides a proper coloring of $\mathbb{T}^d_n$ extending $\xi$; its restriction to the subtree below $v$ is a proper extension of $\xi^v$, again giving $Z_v>0$. At a leaf $v$,
\[
Z_v=1,
\qquad
x_i^{(v)}=\delta_{i,\xi(v)}.
\]

\begin{lemma}[One-step iteration]\label{lem:recursion}
Let $v\in\mathbb{T}^d_n$ satisfy $|v|<n$, so that its children are $v1,\ldots,vd$. Then
\begin{equation}\label{eq:recursion}
x_i^{(v)}
=
\frac{
\displaystyle\prod_{u=1}^d
\bigl[1+(p-1)x_i^{(vu)}\bigr]
}{
\displaystyle\sum_{j=1}^q
\prod_{u=1}^d
\bigl[1+(p-1)x_j^{(vu)}\bigr]
},
\qquad i\in[q].
\end{equation}
\end{lemma}

\begin{proof}
Put $h=n-|v|\geqslant 1$. Fixing the root spin and summing over the configurations of the child subtrees gives
\begin{equation}\label{eq:child-factor}
\begin{aligned}
Z_v x_i^{(v)}
&=
\sum_{\substack{
\sigma\in[q]^{\mathbb{T}^d_h}\\
\sigma_\ast=i
}}
W^{\xi^v}_{h,p,q}(\sigma)\\
&=
\prod_{u=1}^d
\left(
Z_{vu}\sum_{j=1}^q
p^{\delta_{i,j}}x_j^{(vu)}
\right)\\
&=
\left(\prod_{u=1}^d Z_{vu}\right)
\prod_{u=1}^d
\bigl[1+(p-1)x_i^{(vu)}\bigr].
\end{aligned}
\end{equation}
The boundary restriction is included in the weight
$W^{\xi^v}_{h,p,q}$. The last equality uses
\[
\sum_{j=1}^q p^{\delta_{i,j}}x_j^{(vu)}
=
p x_i^{(vu)}+\sum_{j\neq i}x_j^{(vu)}
=
1+(p-1)x_i^{(vu)}.
\]
Summing \eqref{eq:child-factor} over $i$ yields
\[
\sum_{i=1}^q
\prod_{u=1}^d
\bigl[1+(p-1)x_i^{(vu)}\bigr]
=
\frac{Z_v}{\prod_{u=1}^d Z_{vu}}
>0.
\]
Dividing by this sum proves \eqref{eq:recursion}, including at $p=0$.
\end{proof}

\subsection{Centered logarithmic coordinates}

Let
\[
\mathcal{S}_q^\circ
:=
\left\{x\in(0,1)^q:\sum_{i=1}^q x_i=1\right\},
\]
and recall $E_q$ and $\Phi_c$ from \eqref{eq:Eq} and \eqref{eq:Phic}.

The maps
\begin{equation}\label{eq:coords}
y_i(x):=\log x_i-\frac1q\sum_{j=1}^q\log x_j,
\qquad
x_i(y):=\frac{e^{y_i}}{\sum_{j=1}^q e^{y_j}}
\end{equation}
are continuous inverses between $\mathcal{S}_q^\circ$ and $E_q$.
The uniform vector $\pi=(1/q,\ldots,1/q)$ corresponds to $y=0$.

For $p\in[0,1]$, define $T_{d,p}:E_q\to E_q$ by
\begin{equation}\label{eq:T}
T_{d,p,i}(y)
:=
d\left(
\log\bigl(1+(p-1)x_i(y)\bigr)
-\frac1q\sum_{j=1}^q
\log\bigl(1+(p-1)x_j(y)\bigr)
\right),
\qquad i\in[q].
\end{equation}
Since $0<x_i(y)<1$, the logarithms are well defined for every
$p\in[0,1]$, including $p=0$.
The map $T_{d,p}$ is continuous and satisfies $T_{d,p}(0)=0$.

\begin{lemma}[Recursion in centered coordinates]\label{lem:logrec}
Let $v$ be a nonleaf vertex. Suppose that the child root marginals
$x^{(v1)},\ldots,x^{(vd)}$ are strictly positive, and write
$y^{(vu)}=y(x^{(vu)})$. Then $x^{(v)}$ is strictly positive and
\begin{equation}\label{eq:avg}
y^{(v)}
=
\frac1d\sum_{u=1}^d T_{d,p}\bigl(y^{(vu)}\bigr).
\end{equation}
\end{lemma}

\begin{proof}
Every factor $1+(p-1)x_i^{(vu)}$ is positive, since
$x_i^{(vu)}<1$ and $p\in[0,1]$.
Thus \eqref{eq:recursion} gives $x_i^{(v)}>0$ for every $i$.

Let $D_v$ denote the denominator in \eqref{eq:recursion}.
Taking logarithms gives
\[
\log x_i^{(v)}
=
\sum_{u=1}^d
\log\bigl(1+(p-1)x_i^{(vu)}\bigr)
-\log D_v.
\]
The term $\log D_v$ is independent of $i$, so subtracting the
average over colors gives
\[
y_i^{(v)}
=
\sum_{u=1}^d
\left(
\log\bigl(1+(p-1)x_i^{(vu)}\bigr)
-\frac1q\sum_{j=1}^q
\log\bigl(1+(p-1)x_j^{(vu)}\bigr)
\right).
\]
By the definition of $T_{d,p}$, this is \eqref{eq:avg}.
\end{proof}

The child root marginals in \eqref{eq:avg} may differ.
The averaging formula will allow us to apply convexity to
subtrees with arbitrary boundary conditions.
To start this recursion, we need strictly positive marginals
and a common bound on their potentials.
At positive temperature, both follow from the next lemma.

\begin{lemma}[A uniform bound at positive temperature]
\label{lem:initial}
Let $p\in(0,1]$ and let $v$ be a nonleaf vertex.
Then $x^{(v)}$ is strictly positive, even if some coordinates
of the child marginals vanish. Moreover,
\begin{equation}\label{eq:initial-bounds}
p^d
\leqslant
\frac{x_i^{(v)}}{x_j^{(v)}}
\leqslant p^{-d},
\qquad i,j\in[q],
\qquad
\|y^{(v)}\|_\infty
\leqslant d\log\frac1p.
\end{equation}
Consequently, for every $c>0$,
\begin{equation}\label{eq:initial-potential}
\Phi_c\bigl(y^{(v)}\bigr)
\leqslant q\bigl(p^{-cd}-1\bigr).
\end{equation}
\end{lemma}

\begin{proof}
For every $s\in[0,1]$,
\[
p\leqslant1+(p-1)s\leqslant1.
\]
Hence all factors in \eqref{eq:recursion} are positive.
Taking the ratio for colors $i$ and $j$ cancels the common
normalization and gives
\[
\frac{x_i^{(v)}}{x_j^{(v)}}
=
\prod_{u=1}^d
\frac{1+(p-1)x_i^{(vu)}}{1+(p-1)x_j^{(vu)}}
\in[p^d,p^{-d}].
\]
Since
\[
y_i^{(v)}
=
\frac1q\sum_{j=1}^q
\log\frac{x_i^{(v)}}{x_j^{(v)}},
\]
we obtain $|y_i^{(v)}|\leqslant d\log(1/p)$.
Therefore
\[
\Phi_c\bigl(y^{(v)}\bigr)
=
\sum_{i=1}^q e^{cy_i^{(v)}}-q
\leqslant qp^{-cd}-q.
\]
\end{proof}

Thus, when $p>0$, we can start the logarithmic recursion one
level above the leaves, with a potential bound independent
of the boundary condition.

\subsection{A compactness principle}

We first show that a bound on the potential bounds every
coordinate of $y$.

\begin{lemma}[Elementary properties of the potential]
\label{lem:potential}
Let $c>0$. The function $\Phi_c$ is continuous and strictly
convex on $E_q$. It is nonnegative, with equality only at $0$,
and satisfies
\begin{equation}\label{eq:potential-bound}
\|y\|_\infty
\leqslant
\frac{q-1}{c}\sqrt{2\Phi_c(y)},
\qquad y\in E_q.
\end{equation}
In particular, for every $a\geqslant0$, the set
\[
K_a:=\{y\in E_q:\Phi_c(y)\leqslant a\}
\]
is compact.
\end{lemma}

\begin{proof}
Continuity and strict convexity follow from those of the
exponential function. Since $\sum_i y_i=0$,
\begin{equation}\label{eq:potential-centered}
\Phi_c(y)
=
\sum_{i=1}^q\bigl(e^{cy_i}-1-cy_i\bigr).
\end{equation}
Each summand is nonnegative and vanishes only when $y_i=0$.
Hence $\Phi_c(y)\geqslant0$, with equality only at $y=0$.

To prove \eqref{eq:potential-bound}, put $M=\max_i y_i\geqslant0$.
For every $i$,
\[
-y_i=\sum_{j\ne i}y_j\leqslant(q-1)M,
\]
so $\|y\|_\infty\leqslant(q-1)M$.
Keeping a summand corresponding to $M$ in
\eqref{eq:potential-centered}, and using
$e^s-1-s\geqslant s^2/2$ for $s\geqslant0$, gives
\[
\Phi_c(y)
\geqslant e^{cM}-1-cM
\geqslant\frac{c^2M^2}{2}
\geqslant
\frac{c^2\|y\|_\infty^2}{2(q-1)^2}.
\]
This proves \eqref{eq:potential-bound}.
Each $K_a$ is closed by continuity and bounded by this
inequality, so it is compact.
\end{proof}

The next lemma uses this compactness to obtain uniform
convergence from a common finite initial bound, allowing
different vectors at the children of each vertex.

\begin{lemma}[Uniform convergence from global strict descent]
\label{lem:compact}
Let $c>0$ and let $T:E_q\to E_q$ be continuous, with $T(0)=0$ and
\begin{equation}\label{eq:strict-descent}
\Phi_c\bigl(T(y)\bigr)<\Phi_c(y),
\qquad y\in E_q\setminus\{0\}.
\end{equation}
Fix $a_0>0$ and $0<\varepsilon\leqslant a_0$, and define
\begin{equation}\label{eq:compact-constants}
\delta_\varepsilon
:=
\min_{\substack{y\in E_q\\
\varepsilon\leqslant\Phi_c(y)\leqslant a_0}}
\left\{\Phi_c(y)-\Phi_c\bigl(T(y)\bigr)\right\}>0,
\qquad
N_\varepsilon
:=
\left\lceil
\frac{a_0-\varepsilon}{\delta_\varepsilon}
\right\rceil.
\end{equation}
For an integer $m\geqslant0$, assign a vector $y^{(v)}\in E_q$
to each leaf $v$ of $\mathbb{T}^d_m$, with
$\Phi_c(y^{(v)})\leqslant a_0$.
Define the vectors at the remaining vertices recursively by
\[
y^{(v)}
=
\frac1d\sum_{u=1}^d T\bigl(y^{(vu)}\bigr).
\]

Then the root satisfies
\[
\Phi_c\bigl(y^{(\ast)}\bigr)
\leqslant
\max\{\varepsilon,a_0-m\delta_\varepsilon\}.
\]
In particular, its potential is at most $\varepsilon$ whenever
$m\geqslant N_\varepsilon$. Consequently, the root vectors
converge to $0$ as $m\to\infty$, uniformly over all leaf
assignments satisfying the initial bound.
The corresponding probability vectors also converge uniformly
to $(1/q,\ldots,1/q)$.
\end{lemma}

\begin{proof}
\emph{Step 1: a uniform decrease away from $0$.}
The set
\[
A_\varepsilon
:=
\{y\in E_q:
\varepsilon\leqslant\Phi_c(y)\leqslant a_0\}
\]
is compact by Lemma~\ref{lem:potential}.
It is nonempty because
\[
\Phi_c(x,-x,0,\ldots,0)=e^{cx}+e^{-cx}-2
\]
takes every value in $[0,\infty)$ as $x$ runs through
$[0,\infty)$.
Define
\[
g:A_\varepsilon\to\mathbb{R},
\qquad
g(y):=\Phi_c(y)-\Phi_c\bigl(T(y)\bigr).
\]
Since $\Phi_c$ and $T$ are continuous, $g$ is continuous.
As $A_\varepsilon$ is nonempty and compact, $g$ attains
its minimum on $A_\varepsilon$.
Moreover, $0\notin A_\varepsilon$, so
\eqref{eq:strict-descent} gives $g(y)>0$ for every
$y\in A_\varepsilon$. Hence
\[
\delta_\varepsilon
=
\min_{y\in A_\varepsilon}g(y)
>0.
\]
Both $\delta_\varepsilon$ and $N_\varepsilon$ are
independent of the tree height and the leaf assignments.

\emph{Step 2: one iteration with different child vectors.}
Suppose
\[
\Phi_c\bigl(y^{(u)}\bigr)\leqslant a,
\qquad u\in[d],
\qquad
\varepsilon\leqslant a\leqslant a_0.
\]
By \eqref{eq:strict-descent} and $T(0)=0$,
$\Phi_c(T(y))\leqslant\Phi_c(y)$ for every $y\in E_q$.

If $\Phi_c(y^{(u)})<\varepsilon$, then
\[
\Phi_c\bigl(T(y^{(u)})\bigr)
\leqslant\Phi_c(y^{(u)})<\varepsilon.
\]
Otherwise, $y^{(u)}\in A_\varepsilon$, so
\[
\Phi_c\bigl(T(y^{(u)})\bigr)
\leqslant
\Phi_c(y^{(u)})-\delta_\varepsilon
\leqslant a-\delta_\varepsilon.
\]
By convexity,
\begin{equation}\label{eq:compact-step}
\Phi_c\left(\frac1d\sum_{u=1}^d T(y^{(u)})\right)
\leqslant
\frac1d\sum_{u=1}^d\Phi_c\bigl(T(y^{(u)})\bigr)
\leqslant
\max\{\varepsilon,a-\delta_\varepsilon\}.
\end{equation}

\emph{Step 3: iteration and uniform convergence.}
Starting from the initial bound $a_0$, define
\[
a_{j+1}:=\max\{\varepsilon,a_j-\delta_\varepsilon\},
\qquad j\geqslant0.
\]
Applying \eqref{eq:compact-step} successively from the leaves
towards the root shows that every vertex $j$ levels above the
leaves has potential at most $a_j$, for $0\leqslant j\leqslant m$.
The recurrence gives
\[
a_j=\max\{\varepsilon,a_0-j\delta_\varepsilon\}.
\]
Taking $j=m$ proves the stated root bound.
This holds for every height $m$, and the bound is uniform over
all leaf assignments satisfying the initial bound.

In particular,
\[
m\geqslant N_\varepsilon
\quad\Longrightarrow\quad
a_0-m\delta_\varepsilon\leqslant\varepsilon
\quad\Longrightarrow\quad
\Phi_c\bigl(y^{(\ast)}\bigr)\leqslant\varepsilon.
\]

Finally, whenever $\Phi_c(y)\leqslant\varepsilon$,
\eqref{eq:potential-bound} gives
\[
\|y\|_\infty\leqslant\eta_\varepsilon,
\qquad
\eta_\varepsilon:=\frac{q-1}{c}\sqrt{2\varepsilon}.
\]
By \eqref{eq:coords},
\[
\frac{e^{-2\eta_\varepsilon}}q
\leqslant x_i(y)\leqslant
\frac{e^{2\eta_\varepsilon}}q,
\qquad i\in[q].
\]
Consequently,
\begin{equation}\label{eq:x-from-y}
\max_{i\in[q]}\left|x_i(y)-\frac1q\right|
\leqslant
\frac{e^{2\eta_\varepsilon}-1}{q}
=
\frac{\exp\!\left(\frac{2(q-1)}c
\sqrt{2\varepsilon}\right)-1}{q}.
\end{equation}
Both $\eta_\varepsilon$ and the right-hand side tend to
zero as $\varepsilon\downarrow0$. These bounds apply to
the root for every $m\geqslant N_\varepsilon$ and every
leaf assignment satisfying the initial bound.
Since $N_\varepsilon$ is independent of the assignments,
both the root vectors and their corresponding probability
vectors converge uniformly.
\end{proof}

%% file: sections/03-logarithmic-comparison.tex
\section{Logarithmic comparison and the scalar reduction}\label{sec:entropy}

In this section, we reduce the strict decrease of the potential at criticality to an inequality involving a function of one variable. Throughout, fix integers $d\geqslant 2$ and $2\leqslant q\leqslant d+1$. We retain the chosen potential exponent $c>1$, so that $cd>1$. The argument in this section does not require the particular value of $c$ used in Section~\ref{sec:scalar}.

\subsection{Exact identities}

For $y\in E_q$, let $x=x(y)\in\mathcal S_q^\circ$. Since $0<qx_i<q\leqslant d+1$,
\begin{equation}\label{eq:tz}
qx_i>0,\qquad
\frac{d+1-qx_i}{d}>0,
\qquad i\in[q].
\end{equation}
Moreover,
\begin{equation}\label{eq:tz-sums}
\frac{1}{q}\sum_{i=1}^q qx_i
=
\frac{1}{q}\sum_{i=1}^q\frac{d+1-qx_i}{d}
=1.
\end{equation}
Define
\begin{equation}\label{eq:ACBJ}
\begin{aligned}
A&=\frac{1}{q}\sum_{i=1}^q (qx_i)^{c},
&
C&=\frac{1}{q}\sum_{i=1}^q\left(\frac{d+1-qx_i}{d}\right)^{cd},\\
B&=-\frac{c}{q}\sum_{i=1}^q\log(qx_i),
&
J&=-\frac{cd}{q}\sum_{i=1}^q\log\left(\frac{d+1-qx_i}{d}\right).
\end{aligned}
\end{equation}

These quantities express the potential before and after applying
$T_{d,\,1-q/(d+1)}$.

\begin{lemma}\label{lem:identities}
For every $y\in E_q$,
\begin{equation}\label{eq:identities}
\Phi_c(y)+q=qAe^B,
\qquad
\Phi_c\bigl(T_{d,\,1-q/(d+1)}(y)\bigr)+q=qCe^J.
\end{equation}
Moreover, $A,C\geqslant1$ and $B,J\geqslant0$.
\end{lemma}

\begin{proof}
The factor $q$ cancels when we center the logarithms, so
\eqref{eq:coords} gives
\[
y_i
=
\log(qx_i)-\frac1q\sum_{j=1}^q\log(qx_j).
\]
By the definition of $B$,
\[
e^{cy_i}=(qx_i)^c e^B.
\]
Summing over $i$ proves the first identity.

At $p=1-q/(d+1)$, the factors in \eqref{eq:T} satisfy
\[
1+(p-1)x_i
=
1-\frac{qx_i}{d+1}
=
\frac{d}{d+1}\left(\frac{d+1-qx_i}{d}\right).
\]

The common factor $d/(d+1)$ cancels when the logarithms are centered in \eqref{eq:T}. Consequently,
\[
\begin{aligned}
T_{d,1-\frac{q}{d+1},i}(y)
&=
d\left(\log\left(\frac{d+1-qx_i}{d}\right)-\frac{1}{q}\sum_{j=1}^q\log\left(\frac{d+1-qx_j}{d}\right)\right),\\
e^{cT_{d,1-\frac{q}{d+1},i}(y)}&=\left(\frac{d+1-qx_i}{d}\right)^{cd}e^J.
\end{aligned}
\]
Summing over $i$ gives the second identity.

Since $c,cd>1$, the functions $s\mapsto s^{c}$ and
$s\mapsto s^{cd}$ are convex on $(0,\infty)$, as is
$s\mapsto-\log s$. Hence, by \eqref{eq:tz-sums},
Jensen's inequality gives
\[
\begin{aligned}
A&=\frac1q\sum_{i=1}^q(qx_i)^{c}
\ge \left(\frac1q\sum_{i=1}^q qx_i\right)^{c}=1,\\
C&=\frac1q\sum_{i=1}^q\left(\frac{d+1-qx_i}{d}\right)^{cd}
\ge \left(\frac1q\sum_{i=1}^q\frac{d+1-qx_i}{d}\right)^{cd}=1,
\end{aligned}
\]
and
\[
\begin{aligned}
B&=-\frac{c}{q}\sum_{i=1}^q\log(qx_i)
\ge -c\log\left(\frac1q\sum_{i=1}^q qx_i\right)=0,\\
J&=-\frac{cd}{q}\sum_{i=1}^q\log\left(\frac{d+1-qx_i}{d}\right)\\
&\ge -cd\log\left(\frac1q\sum_{i=1}^q\frac{d+1-qx_i}{d}\right)=0.
\end{aligned}
\]
\end{proof}

\subsection{Comparison of the logarithmic terms}

We now prove $J\leqslant B$ directly. By \eqref{eq:ACBJ},
\[
J-B
=
\frac cq\sum_{i=1}^q
\left[
\log(qx_i)
-d\log\left(\frac{d+1-qx_i}{d}\right)
\right].
\]
We will show that this expression is maximized at the uniform vector,
where its value is zero.

\begin{lemma}[Logarithmic comparison]\label{lem:Kprops}
For every $x\in\mathcal S_q^\circ$,
\begin{equation}\label{eq:JleB}
J\leqslant B.
\end{equation}
\end{lemma}

\begin{proof}
We prove that $J-B$ is maximized at the uniform vector, where
its value is zero. For this argument, allow $d$ to be real,
keeping $q$ and $c$ fixed, and first assume that $d>1$ and
$q<d+1$.

Recall that
\[
J-B
=
\frac cq\sum_{i=1}^q
\left[
\log(qx_i)
-d\log\left(\frac{d+1-qx_i}{d}\right)
\right].
\]
Since $0<x_i<1$,
\[
0<
\frac{d+1-q}{d}
\leqslant
\frac{d+1-qx_i}{d}
\leqslant
\frac{d+1}{d}.
\]
Thus the logarithms involving $(d+1-qx_i)/d$ are uniformly
bounded. On the other hand, $\log(qx_i)\leqslant\log q$,
and the smallest of these terms tends to $-\infty$ as
$\min_i x_i\to0$. Hence
\[
J-B\longrightarrow-\infty
\qquad\text{as }\min_i x_i\longrightarrow0.
\]

At the uniform vector, $J-B=0$. We may therefore choose
$0<\varepsilon<1/q$ such that $J-B<0$ whenever
$\min_i x_i<\varepsilon$.
The continuous function $J-B$ attains a maximum on the
closed and bounded set
\[
\left\{
x\in\mathbb R^q:
\sum_{i=1}^q x_i=1,\quad
x_i\geqslant\varepsilon\ \text{for every }i
\right\}.
\]
This maximum is at least zero, whereas all points outside
this set have negative values. It is therefore a global
maximum, attained at a point with all coordinates positive.

We now show that such a maximizing point must be uniform.
The coordinate derivatives are
\[
\begin{aligned}
\frac{\partial(J-B)}{\partial x_i}
&=
c\left(
\frac1{qx_i}+\frac d{d+1-qx_i}
\right),\\
\frac{\partial^2(J-B)}{\partial x_i^2}
&=
cq\left(
-\frac1{(qx_i)^2}
+\frac d{(d+1-qx_i)^2}
\right).
\end{aligned}
\]
At a maximizing point, increasing $x_i$ and decreasing $x_j$
by the same small amount preserves the constraint
$\sum_i x_i=1$. The derivative along this variation must
vanish, so
\[
\frac1{qx_i}+\frac d{d+1-qx_i}
=
\frac1{qx_j}+\frac d{d+1-qx_j}
\qquad(i,j\in[q]).
\]

On the interval $0<x_i<(d+1)/q$, the second partial derivative
is negative precisely when
\[
\sqrt d\,qx_i<d+1-qx_i,
\qquad\text{that is,}\qquad
x_i<\frac{d+1}{q(1+\sqrt d)},
\]
and positive when the inequality is reversed.
Thus the first partial derivative, viewed as a function of
its coordinate, is strictly decreasing below this point and
strictly increasing above it.

Thus each value of the first derivative is attained at most once
on either side of the dividing point. At an interior maximum
under the constraint $\sum_{i=1}^q x_i=1$, all first partial
derivatives are equal. Hence the coordinates take at most two
distinct values. If there are two, one lies strictly below the
dividing point and the other strictly above it.

Suppose that two distinct coordinate values occur. We claim
that the larger value occurs only once. Otherwise, choose two
coordinates $x_i$ and $x_j$ with this value and consider the
variation
\[
x_i\mapsto x_i+h,
\qquad
x_j\mapsto x_j-h.
\]
For sufficiently small $|h|$, the coordinates remain positive
and their sum remains one. Since $J-B$ is a sum of functions
of individual coordinates, its mixed second derivatives vanish.
The second derivative along this variation at $h=0$ is therefore
\[
\frac{\partial^2(J-B)}{\partial x_i^2}
+
\frac{\partial^2(J-B)}{\partial x_j^2}
>0,
\]
because both coordinates lie above the dividing point.
This contradicts the second derivative condition at a maximum.

Consequently, after relabelling, any nonuniform maximizing point
must satisfy
\[
x_1=\cdots=x_{q-1}<\frac1q<x_q,
\qquad
x_q=1-(q-1)x_1.
\]
Equality of the first partial derivatives corresponding to
$x_1$ and $x_q$ gives
\[
0
=
q(x_q-x_1)
\left[
\frac{1}{q^2x_1x_q}
-
\frac{d}{(d+1-qx_1)(d+1-qx_q)}
\right].
\]
Since $x_q>x_1$ and all denominators are positive, it follows that
\[
(d+1-qx_1)(d+1-qx_q)-dq^2x_1x_q=0.
\]
On the other hand, substituting $x_q=1-(q-1)x_1$ gives
\begin{align*}
&(d+1-qx_1)(d+1-qx_q)-dq^2x_1x_q\\
&\qquad
=(d+1-q)(d+1-2qx_1)
 +(q-1)(d-1)q^2x_1^2
>0.
\end{align*}
Indeed, $d+1-q\geqslant0$ and $qx_1<1$ imply that
the first term is nonnegative, since
\[
d+1-2qx_1>d-1>0.
\]
The second term is strictly positive.

This contradiction rules out every nonuniform maximizing point.
The maximum is therefore attained at the uniform vector,
where $J-B=0$. Hence $J\leqslant B$ when $q<d+1$.

Finally, suppose that $q=d+1$.
Fix $q$, $c$, and a positive probability vector $x$, and let
the real parameter $d$ decrease to $q-1$ from above.
The inequality just proved holds before taking the limit.
At the limiting value,
\[
\frac{d+1-qx_i}{d}
=
\frac{q(1-x_i)}{q-1}>0.
\]
Thus all logarithms remain continuous at the limit, and
passing to the limit in $J-B\leqslant0$ proves
\eqref{eq:JleB} also when $q=d+1$.
\end{proof}

\subsection{Reduction to a scalar inequality}

The estimate $J\leqslant B$ controls the exponential factors in \eqref{eq:identities}. We must also compare $A$ and $C$.

For $0<x<(d+1)/q$, define
\begin{equation}\label{eq:Fdef}
F_{q,d}(x)
:=
(qx)^c
-\left(\frac{d+1-qx}{d}\right)^{cd}
-c\log(qx)
+cd\log\left(\frac{d+1-qx}{d}\right).
\end{equation}
The definitions of $A,C,B,J$ give
\begin{equation}\label{eq:F-sum}
\frac1q\sum_{i=1}^q F_{q,d}(x_i)
=
A-C+B-J.
\end{equation}
The next lemma shows that positivity of this sum implies
strict decrease of the potential.

\begin{lemma}[Scalar sufficient condition]\label{lem:scalar}
Let $y\in E_q$ and write $x=x(y)$. If
\[
\sum_{i=1}^q F_{q,d}(x_i)>0,
\]
then
\[
\Phi_c\bigl(T_{d,\,1-q/(d+1)}(y)\bigr)<\Phi_c(y).
\]
\end{lemma}

\begin{proof}
By \eqref{eq:identities}, and since $A,C>0$,
\begin{align*}
\Phi_c\bigl(T_{d,\,1-q/(d+1)}(y)\bigr)<\Phi_c(y)
&\quad\Longleftrightarrow\quad Ce^J<Ae^B\\
&\quad\Longleftrightarrow\quad
\log\frac{C}{A}<B-J.
\end{align*}
It therefore suffices to prove the last inequality.
The hypothesis and \eqref{eq:F-sum} give
\[
C-A<B-J,
\]
while the preceding lemma gives $B-J\geqslant0$.
If $C>A$, then $A\geqslant 1$ and $\log(1+u)\leqslant u$ imply
\[
\log\frac{C}{A}
\leqslant
\frac{C-A}{A}
\leqslant
C-A
<
B-J.
\]
If $C\leqslant A$, then
\[
\log\frac{C}{A}\leqslant 0\leqslant B-J.
\]
Equality between the two ends would require both $C=A$ and $B=J$, contradicting $A-C+B-J>0$. The desired strict inequality therefore holds in both cases.
\end{proof}

Since $y=0$ if and only if $x_1=\cdots=x_q=1/q$, it remains to prove
\[
\sum_{i=1}^q F_{q,d}(x_i)>0
\]
for every positive vector $(x_1,\ldots,x_q)$ with $\sum_i x_i=1$ other than $(1/q,\ldots,1/q)$. This is the main task of Section~\ref{sec:scalar}.

%% file: sections/04-scalar.tex
\section{The scalar inequality and critical descent}\label{sec:scalar}

Throughout this section, let $d\geqslant2$ and
$2\leqslant q\leqslant d+1$ be integers.

For the moment, let $c>1$ be arbitrary and write $F=F_{q,d}$
for the function defined in \eqref{eq:Fdef}:
\begin{equation}
F(x)
=
(qx)^c
-\left(\frac{d+1-qx}{d}\right)^{cd}
-c\log(qx)
+cd\log\left(\frac{d+1-qx}{d}\right),
\qquad
0<x<\frac{d+1}{q}.
\end{equation}
Our goal is to prove
\[
\sum_{i=1}^q F_{q,d}(x_i)>0
\]
for every positive probability vector other than
$(1/q,\ldots,1/q)$. We will first establish
\[
F(x)>0
\qquad
\left(0<x\leqslant\frac dq,\quad x\ne\frac1q\right).
\]
Since $d\geqslant2$, this requires positivity on both sides
of $x=1/q$. To choose the exponent $c$, we examine the
Taylor expansion at this point.

Using
\[
\frac{d}{dx}\left(\frac{d+1-qx}{d}\right)=-\frac qd,
\]
we obtain
\[
F'(x)
=
cq\left[
(qx)^{c-1}
+\left(\frac{d+1-qx}{d}\right)^{cd-1}
-\frac1{qx}
-\frac d{d+1-qx}
\right].
\]
At $x=1/q$, both $qx$ and $(d+1-qx)/d$ equal $1$, so
\[
F(1/q)=0,
\qquad
F'(1/q)=cq(1+1-1-1)=0.
\]
Differentiating further gives
\[
F''(1/q)
=
cq^2\left[
(c-1)-\frac{cd-1}{d}+1-\frac1d
\right]
=0
\]
and
\begin{align*}
F'''(1/q)
&=
cq^3\left[
(c-1)(c-2)
+\frac{(cd-1)(cd-2)}{d^2}
-2-\frac2{d^2}
\right]\\
&=
2c^2q^3\left(c-\frac{3(d+1)}{2d}\right).
\end{align*}
Thus, as $h\to0$,
\begin{align*}
F(1/q+h)
&=
F(1/q)+F'(1/q)h
+\frac{F''(1/q)}2h^2
+\frac{F'''(1/q)}6h^3
+O(h^4)\\
&=
\frac{c^2q^3}{3}
\left(c-\frac{3(d+1)}{2d}\right)h^3
+O(h^4).
\end{align*}
If the cubic coefficient is nonzero, it determines the sign
of $F(1/q+h)$ for sufficiently small nonzero $h$.
Choosing $h$ with the opposite sign to this coefficient
would give $F(1/q+h)<0$. Since both signs of $h$ are
allowed near $1/q$, the pointwise positivity we seek
requires $F'''(1/q)=0$. As $c>1$, this determines
\begin{equation}
c=\frac{3(d+1)}{2d},
\qquad
cd=\frac{3(d+1)}2.
\end{equation}
We fix this value of $c$ for the rest of the section.
In particular,
\[
\frac32<c\leqslant\frac94,
\qquad
cd\geqslant\frac92.
\]

This cancellation explains the choice of $c$, but the
required positivity still needs to be proved. Moreover,
\[
F(x)\longrightarrow-\infty
\qquad\text{as }x\uparrow\frac{d+1}{q},
\]
because
\[
cd\log\left(\frac{d+1-qx}{d}\right)\longrightarrow-\infty,
\]
while the other terms in the definition of $F$ remain bounded.
We must therefore use the constraint $\sum_i x_i=1$.

We first prove that $F_{q,d}(x)>0$ for
$0<x\leqslant d/q$, except at $x=1/q$.
This settles the case $q\leqslant d$, since every coordinate
then satisfies $x_i<1\leqslant d/q$.
When $q=d+1$, at most one coordinate can exceed $d/q$.
In that case, we use the remaining coordinates to compensate
for its possible negative contribution.

\subsection{Shape of the scalar function}

We first determine where $F'$ changes sign and prove that
$F$ is convex on $(0,1/q)$. This convexity will later allow
us to apply Jensen's inequality to the smaller coordinates.
In the proof, we introduce an auxiliary function $P$ by
removing the denominators from $F'$ without changing its sign.
The expression for $P'''$ consists of two powers multiplied
by affine functions, whose signs and magnitudes we can compare.

\begin{lemma}[Shape of $F_{q,d}$]\label{lem:shape}
For all integers $d\geqslant2$ and $2\leqslant q\leqslant d+1$,
the following statements hold.
\begin{enumerate}
\item[(i)]
For $0<x<1/q$,
\[
F'_{q,d}(x)<0,
\qquad
F''_{q,d}(x)>0,
\qquad
F_{q,d}(x)>0.
\]

\item[(ii)]
There is a unique $\rho\in(1/q,(d+1)/q)$ such that
\[
F'_{q,d}(x)>0
\qquad (1/q<x<\rho),
\]
and
\[
F'_{q,d}(x)<0
\qquad (\rho<x<(d+1)/q).
\]

\item[(iii)]
If $F_{q,d}(d/q)>0$, then
\[
F_{q,d}(x)>0
\qquad
\left(0<x\leqslant\frac dq,\quad x\ne\frac1q\right).
\]
\end{enumerate}
\end{lemma}

\begin{proof}
Differentiating the definition of $F_{q,d}$ gives
\begin{equation}\label{eq:F1}
F'_{q,d}(x)
=
\frac{3q(d+1)}{2d}
\left[
(qx)^{\frac{d+3}{2d}}
+\left(\frac{d+1-qx}{d}\right)^{\frac{3d+1}{2}}
-\frac1{qx}
-\frac d{d+1-qx}
\right].
\end{equation}
At $x=1/q$, further differentiation gives
\begin{equation}\label{eq:F4}
F_{q,d}(1/q)
=
F'_{q,d}(1/q)
=
F''_{q,d}(1/q)
=
F'''_{q,d}(1/q)
=0,
\end{equation}
and
\[
F^{(4)}_{q,d}(1/q)
=
\frac{9q^4(d+1)^2(d^2-1)}{2d^4}
>0.
\]

To determine the sign of $F_{q,d}'$, define
\begin{equation}\label{eq:Pdef}
\begin{aligned}
P(x)
&=
\frac{2x(d+1-qx)}{3(d+1)}F_{q,d}'(x)\\
&=
\frac{d+1-qx}{d}\,(qx)^{\frac{3(d+1)}{2d}}
+
qx\left(\frac{d+1-qx}{d}\right)^{\frac{3(d+1)}2}
-qx-\frac{d+1-qx}{d}.
\end{aligned}
\end{equation}
The factor multiplying $F_{q,d}'$ is positive on $(0,(d+1)/q)$,
so $P$ and $F_{q,d}'$ have the same sign.
The product rule and \eqref{eq:F4} give
\begin{equation}\label{eq:P-at-1}
P(1/q)=P'(1/q)=P''(1/q)=0.
\end{equation}

Differentiating \eqref{eq:Pdef} three times gives
\begin{equation}\label{eq:P3}
\begin{aligned}
P'''(x)
={}&
-\frac{3q^3(d+1)(d+3)}{8d^4}\,
(qx)^{-\frac{3(d-1)}{2d}}
\bigl[(d-3)(d+1)+(5d+3)qx\bigr]\\
&+
\frac{3q^3(d+1)(3d+1)}{8d^3}
\left(\frac{d+1-qx}{d}\right)^{\frac{3(d-1)}2}
\bigl[6(d+1)-(3d+5)qx\bigr].
\end{aligned}
\end{equation}
We will also use the endpoint limits
\begin{equation}\label{eq:P-endpoints}
\begin{aligned}
P(((d+1)/q)^-)
&=-(d+1)<0,\\
P'(((d+1)/q)^-)
&=q\left[-\frac{(d+1)^{\frac{3(d+1)}{2d}}}{d}
-1+\frac1d\right]<0,\\
P''(((d+1)/q)^-)
&=-\frac{3q^2}{d^2}(d+1)^{\frac{3(d+1)}{2d}}<0,\\
P'(0^+)
&=q\left[\left(1+\frac1d\right)^{\frac{3(d+1)}2}
-1+\frac1d\right]>0.
\end{aligned}
\end{equation}
For the last inequality, the base is greater than $1$ and
the exponent is positive, so the power is greater than $1$.

We first prove that $P'''$ is positive on a single interval
containing $1/q$ and negative outside its closure.

\emph{The case $d\geqslant3$.}
In \eqref{eq:P3},
\[
(d-3)(d+1)+(5d+3)qx>0
\qquad(x>0),
\]
so the first term is strictly negative.
When $x\geqslant6(d+1)/(q(3d+5))$, the bracket in the second
term is nonpositive. Hence
\[
P'''(x)<0
\qquad
\left(\frac{6(d+1)}{q(3d+5)}\leqslant x<(d+1)/q\right).
\]

For $0<x<6(d+1)/(q(3d+5))$, the second term is positive.
Its ratio to the magnitude of the first term is
\begin{equation}\label{eq:R}
R(x)
=
\frac{d(3d+1)}{d+3}\,
(qx)^{\frac{3(d-1)}{2d}}
\left(\frac{d+1-qx}{d}\right)^{\frac{3(d-1)}2}
\frac{6(d+1)-(3d+5)qx}
     {(d-3)(d+1)+(5d+3)qx}.
\end{equation}

All factors in this expression are positive on the interval
under consideration. By definition, $R$ is the ratio of the
positive term in $P'''$ to the absolute value of the negative
term. Hence $P'''$ has the same sign as $R-1$.

Differentiating $\log R$ twice gives
\begin{align*}
(\log R)''(x)
={}&-\frac{3(d-1)}{2dx^2}
-\frac{3q^2(d-1)}{2(d+1-qx)^2}\\
&-\frac{q^2(3d+5)^2}
       {[6(d+1)-(3d+5)qx]^2}
+\frac{q^2(5d+3)^2}
       {[(d-3)(d+1)+(5d+3)qx]^2}.
\end{align*}
To bound the positive term, note that $d\geqslant3$ gives
\[
(d-3)(d+1)+(5d+3)qx
\geqslant(5d+3)qx>0.
\]
It follows that
\[
\frac{q^2(5d+3)^2}
     {[(d-3)(d+1)+(5d+3)qx]^2}
\leqslant\frac1{x^2}.
\]
Consequently,
\begin{equation}
(\log R)''(x)
\leqslant
-\frac{d-3}{2dx^2}
-\frac{3q^2(d-1)}{2(d+1-qx)^2}
-\frac{q^2(3d+5)^2}
       {[6(d+1)-(3d+5)qx]^2}
<0.
\end{equation}
The first term is nonpositive, and the other two are strictly
negative. Thus $\log R$ is strictly concave. We also have
\[
R(1/q)=\left(\frac{3d+1}{d+3}\right)^2>1,
\qquad
\lim_{x\uparrow 6(d+1)/(q(3d+5))}R(x)=0.
\]

By concavity, the set where $\log R>0$ is an interval:
between any two points in this set, the graph lies above
the positive chord joining them. This interval contains
$1/q$. Its right endpoint is strictly below
$6(d+1)/(q(3d+5))$, since $\log R$ tends to $-\infty$
at that point.

Outside the closure of the positive interval, $\log R$ is
strictly negative. Otherwise, it would vanish at some such
point. Strict concavity between that point and $1/q$ would
then give positive values outside the closure, a contradiction.
Since $P'''$ has the same sign as $R-1$ and is already known
to be negative to the right of $6(d+1)/(q(3d+5))$, this proves
the asserted sign pattern.

\emph{The case $d=2$.}
The two brackets in the expression for $P'''$ become
$13qx-3$ and $18-11qx$, respectively.
With the minus sign before the first term, this gives
\[
P'''>0
\quad\text{on }\left(0,\frac{3}{13q}\right],
\qquad
P'''<0
\quad\text{on }\left[\frac{18}{11q},\frac3q\right).
\]
On $(3/(13q),18/(11q))$, the first term is negative and
the second is positive. The ratio of the positive term
to the absolute value of the negative term is
\[
R(x)
=
\frac{14}{5}(qx)^{3/4}
\left(\frac{3-qx}{2}\right)^{3/2}
\frac{18-11qx}{13qx-3}.
\]
Its logarithmic derivative is
\[
(\log R)'(x)
=
\frac{3}{4x}
-\frac{3q}{2(3-qx)}
-\frac{11q}{18-11qx}
-\frac{13q}{13qx-3}.
\]
All denominators are positive on this interval. Moreover,
\[
13qx-3<13qx
\quad\Longrightarrow\quad
\frac{13q}{13qx-3}>\frac1x>\frac{3}{4x}.
\]
Thus the sum of the first and last terms is negative.
The other two terms are also negative, so $(\log R)'<0$.

Since
\[
\lim_{x\downarrow 3/(13q)}R(x)=+\infty,
\qquad
R(1/q)=\frac{49}{25}>1,
\qquad
\lim_{x\uparrow 18/(11q)}R(x)=0,
\]
the function $R$ crosses $1$ exactly once, at a point
greater than $1/q$. Hence $P'''$ has the asserted sign
pattern also when $d=2$.

We now use this sign pattern to determine the signs of
$P$ and $F'_{q,d}$.

\emph{The interval $(0,1/q)$.}
For $x<1/q$ within the interval where $P'''>0$,
\eqref{eq:P-at-1} gives
\[
P''(x)=-\int_x^{1/q}P'''(s)\,ds<0,
\qquad
P'(x)=-\int_x^{1/q}P''(s)\,ds>0.
\]
If this positive interval has a left endpoint greater than
$0$, the same integrals show that $P'$ is strictly positive
at that endpoint.

On the remaining interval between $0$ and that endpoint,
$P'''<0$.
Since $(P')''=P'''$, the function $P'$ is concave there.
Its two endpoint values are positive: the left one is
$P'(0^+)>0$ by \eqref{eq:P-endpoints}, and the right one
was obtained by integration.
The graph of $P'$ lies above the chord joining these values,
so $P'>0$ on this remaining interval as well.

Hence $P'(x)>0$ for every $0<x<1/q$. Since $P(1/q)=0$,
\[
P(x)=-\int_x^{1/q}P'(s)\,ds<0.
\]
The sign relation between $P$ and $F'_{q,d}$ then gives
$F'_{q,d}(x)<0$. Using $F_{q,d}(1/q)=0$, we obtain
\[
F_{q,d}(x)
=
-\int_x^{1/q}F'_{q,d}(s)\,ds
>0.
\]

To prove convexity, we write
\[
F'_{q,d}(x)
=
\frac{3(d+1)}2
\frac{P(x)}{x(d+1-qx)}.
\]
Differentiating gives
\begin{equation}
F''_{q,d}(x)
=
\frac{3(d+1)}2
\frac{x(d+1-qx)P'(x)-(d+1-2qx)P(x)}
     {x^2(d+1-qx)^2}.
\end{equation}
For $0<x<1/q$, the denominator is positive, and
\[
x(d+1-qx)>0,
\qquad
d+1-2qx>d-1>0.
\]
Since $P'(x)>0$ and $P(x)<0$, both terms in the numerator
are positive. Thus $F''_{q,d}(x)>0$, proving (i).

\emph{The interval $(1/q,(d+1)/q)$.}
On this interval, $P'''$ changes sign exactly once, from
positive to negative. Hence $P''$ first increases strictly
from $P''(1/q)=0$ and then decreases strictly.
Its maximum is positive, while its limit at the right
endpoint is negative by \eqref{eq:P-endpoints}. It therefore crosses zero
exactly once, from positive to negative.

It follows that $P'$ first increases strictly from
$P'(1/q)=0$ and then decreases strictly.
Again, its maximum is positive and its limit at the right
endpoint is negative, so it crosses zero exactly once,
from positive to negative.

In turn, $P$ first increases strictly from $P(1/q)=0$
and then decreases strictly to its negative limit at the
right endpoint. Thus there is a unique
$\rho\in(1/q,(d+1)/q)$ such that
\[
P>0
\quad\text{on }(1/q,\rho),
\qquad
P<0
\quad\text{on }(\rho,(d+1)/q).
\]
Since $F'_{q,d}$ has the same sign as $P$, this proves (ii).

Finally, suppose that $F_{q,d}(d/q)>0$.
Part (i) gives positivity on $(0,1/q)$.
On $(1/q,\rho]$, part (ii) gives
\[
F_{q,d}(x)>F_{q,d}(1/q)=0.
\]
If $d/q\leqslant\rho$, this proves positivity throughout
$(1/q,d/q]$. If $d/q>\rho$, then $F_{q,d}$ decreases on
$[\rho,d/q]$, so
\[
F_{q,d}(x)\geqslant F_{q,d}(d/q)>0
\qquad
(\rho\leqslant x\leqslant d/q).
\]
This proves (iii).
\end{proof}

\subsection{Endpoint positivity}

By the shape lemma, positivity on
$(0,d/q]\setminus\{1/q\}$ reduces to proving
$F_{q,d}(d/q)>0$. We treat $d=2$ directly.
For $d\geqslant3$, we give a uniform argument by comparing
the positive power term with the logarithmic term.

\begin{lemma}[Endpoint positivity]\label{lem:endpoint}
For every integer $d\geqslant2$,
\begin{equation}\label{eq:endpoint}
F_{q,d}(d/q)>0.
\end{equation}
Consequently,
\begin{equation}\label{eq:F-pos}
F_{q,d}(x)>0
\qquad(0<x\leqslant d/q,\ x\neq1/q).
\end{equation}
\end{lemma}

\begin{proof}
By definition,
\begin{equation}\label{eq:Fdd}
F_{q,d}(d/q)
=
d^{3(d+1)/(2d)}
-d^{-3(d+1)/2}
-\frac{3(d+1)^2}{2d}\log d.
\end{equation}

For $d=2$, the elementary bounds
\[
2^{1/4}>\frac{19}{16},
\qquad
\sqrt2>\frac75,
\qquad
\log2<\frac{25}{36}
\]
give
\[
\begin{aligned}
F_{q,2}(2/q)
&=
4\cdot2^{1/4}
-\frac1{16\sqrt2}
-\frac{27}{4}\log2\\
&>
\frac{19}{4}-\frac5{112}-\frac{75}{16}
>0.
\end{aligned}
\]

For the remaining case, regard the right-hand side of
\eqref{eq:Fdd} as a function of real $d\geqslant3$.
We first show that
\[
\frac{d^{3(d+1)/(2d)}}
     {\frac{3(d+1)^2}{2d}\log d}
\]
is strictly increasing.
This will give a lower bound for the difference between
the first and third terms in \eqref{eq:Fdd}.

We use the logarithmic bounds
\[
\frac{3(d-1)}{d+3}<\log d<\frac{d+1}{3}
\qquad(d\geqslant3).
\]
Indeed,
\[
\frac{\mathrm d}{\mathrm d d}
\left(\log d-\frac{3(d-1)}{d+3}\right)
=
\frac{(d-3)^2}{d(d+3)^2}\geqslant0,
\]
and
\[
\frac{\mathrm d}{\mathrm d d}
\left(\frac{d+1}{3}-\log d\right)
=
\frac13-\frac1d\geqslant0.
\]
Both differences are strictly positive at $d=3$, since
\[
1<\log3<\frac{11}{10}<\frac43.
\]
The stated bounds therefore hold throughout $[3,\infty)$.

Differentiating the logarithm of the ratio gives
\[
\frac{\mathrm d}{\mathrm d d}
\log\left(
\frac{d^{3(d+1)/(2d)}}
     {\frac{3(d+1)^2}{2d}\log d}
\right)
=
\frac{
-3(d+1)(\log d)^2
+(d^2+8d+3)\log d
-2d(d+1)}
{2d^2(d+1)\log d}.
\]
The denominator is positive.
For fixed $d$, the numerator is a strictly concave quadratic
in the quantity occupying the position of $\log d$.
At the two endpoints of the interval
\[
\left[\frac{3(d-1)}{d+3},\,\frac{d+1}{3}\right],
\]
its values are
\[
\begin{aligned}
&-3(d+1)\left(\frac{3(d-1)}{d+3}\right)^2
+(d^2+8d+3)\frac{3(d-1)}{d+3}
-2d(d+1)\\
&\qquad=
\frac{(d-3)\bigl((d-3)^3+(d-3)^2+36\bigr)}
     {(d+3)^2}
\geqslant0
\end{aligned}
\]
and
\[
\begin{aligned}
&-3(d+1)\left(\frac{d+1}{3}\right)^2
+(d^2+8d+3)\frac{d+1}{3}
-2d(d+1)\\
&\qquad=\frac{2(d+1)}3>0.
\end{aligned}
\]
A concave function lies above the chord joining its endpoint
values. Since $\log d$ lies strictly between these two endpoints,
and the second endpoint value is strictly positive,
the numerator is strictly positive.
Thus the ratio is strictly increasing, and
\[
\frac{d^{3(d+1)/(2d)}}
     {\frac{3(d+1)^2}{2d}\log d}
\geqslant
\frac9{8\log3}
>1
\qquad(d\geqslant3).
\]

Moreover,
\[
\frac{(d+1)^2}{d}=d+2+\frac1d
\]
is positive and increasing for $d\geqslant3$, because its derivative
is $1-d^{-2}>0$.
Since $\log d$ is also positive and increasing,
\[
\frac{3(d+1)^2}{2d}\log d\geqslant8\log3.
\]
Both factors in the following product therefore have positive
lower bounds:
\[
\begin{aligned}
&d^{3(d+1)/(2d)}
-\frac{3(d+1)^2}{2d}\log d\\
&\quad=
\frac{3(d+1)^2}{2d}\log d
\left(
\frac{d^{3(d+1)/(2d)}}
     {\frac{3(d+1)^2}{2d}\log d}
-1
\right)\\
&\quad\geqslant
8\log3\left(\frac9{8\log3}-1\right)
=
9-8\log3.
\end{aligned}
\]

Finally, $d\geqslant3$ and $3(d+1)/2\geqslant6$ imply
\[
d^{-3(d+1)/2}\leqslant3^{-6}=\frac1{729}.
\]
Combining these estimates with \eqref{eq:Fdd}, and using
$\log3<11/10$, gives
\[
F_{q,d}(d/q)
\geqslant
9-8\log3-\frac1{729}
>
\frac15-\frac1{729}
>0.
\]
This proves \eqref{eq:endpoint}.
Equation \eqref{eq:F-pos} now follows from
Lemma~\ref{lem:shape}(iii).
\end{proof} 
\paragraph{The case with no large coordinate.}
Let $x=(x_1,\ldots,x_q)$ be a positive, nonuniform probability vector.
If $\max_i x_i\leqslant d/q$, then \eqref{eq:F-pos} and $F_{q,d}(1/q)=0$
give $F_{q,d}(x_i)\geqslant0$ for every $i$.
At least one coordinate differs from $1/q$, so
\[
\sum_{i=1}^q F_{q,d}(x_i)>0.
\]
For $2\leqslant q\leqslant d$, every such vector satisfies
$x_i<1\leqslant d/q$ for all $i$.
Lemma~\ref{lem:scalar} therefore gives strict potential descent at the
positive critical parameter throughout this range.

It remains to consider $q=d+1$ with $\max_i x_i>d/q$.

\input{sections/04-compensation.tex}

\subsection{The constrained sum and critical descent}

\begin{theorem}[Constrained scalar inequality]\label{thm:constrained}
Let $d\geqslant2$ and $2\leqslant q\leqslant d+1$ be integers.
If $x_1,\ldots,x_q>0$,
\[
\sum_{i=1}^q x_i=1,
\qquad
(x_1,\ldots,x_q)\neq(1/q,\ldots,1/q),
\]
then
\begin{equation}\label{eq:constrained}
\sum_{i=1}^q F_{q,d}(x_i)>0.
\end{equation}
\end{theorem}

\begin{proof}
If $\max_i x_i\leqslant d/q$, the claim follows from \eqref{eq:F-pos},
$F_{q,d}(1/q)=0$, and the nonuniformity of $x$.
Otherwise $q=d+1$, and we may relabel the coordinates so that $x_1>d/q$.
The Jensen bound in Section~\ref{subsec:compensation} and
Lemma~\ref{lem:compensation} give
\[
\sum_{i=1}^{d+1}F_{q,d}(x_i)\geqslant S_{q,d}(x_1)>0.
\]
\end{proof}

\begin{theorem}[Strict descent at the critical parameter]\label{thm:descent}
Let $d\geqslant2$ and $2\leqslant q\leqslant d+1$ be integers.
Then
\begin{equation}\label{eq:descent}
\Phi_{\frac{3(d+1)}{2d}}\bigl(T_{d,1-\frac{q}{d+1}}(y)\bigr)
<
\Phi_{\frac{3(d+1)}{2d}}(y)
\qquad
\text{for every }y\in E_q\setminus\{0\}.
\end{equation}
\end{theorem}

\begin{proof}
For $y\in E_q\setminus\{0\}$, let $x=x(y)$.
The coordinates $x_i$ are positive with sum $1$.
They are not all equal to $1/q$, since $x(y)=\pi$ if and only if $y=0$.
Theorem~\ref{thm:constrained} therefore gives
\[
\sum_{i=1}^qF_{q,d}(x_i)>0,
\]
and Lemma~\ref{lem:scalar} yields \eqref{eq:descent}.
\end{proof}

\begin{remark}\label{rem:zero-temp}
When $q=d+1$, the critical parameter is $1-q/(d+1)=0$.
Theorem~\ref{thm:descent} gives strict descent for every strictly positive root marginal other than the uniform vector; at the uniform vector, both potentials are zero.
At zero temperature, however, an admissible boundary condition
may give a root marginal with zero coordinates.
The centered logarithmic coordinates are then not finite,
so the theorem does not establish Gibbs uniqueness at
$(q,p)=(d+1,0)$.
\end{remark}

%% file: sections/04-compensation.tex
\subsection{Compensation when $q=d+1$}\label{subsec:compensation}

Suppose $q=d+1$ and one coordinate exceeds $d/q$.
After relabelling, let $x_1>d/q$.
The remaining $d$ coordinates sum to $1-x_1<1/q$,
so each belongs to $(0,1/q)$.
Since $F_{q,d}$ is convex on this interval by
Lemma~\ref{lem:shape}(i), Jensen's inequality gives
\[
\sum_{i=1}^{d+1}F_{q,d}(x_i)
\geqslant
F_{q,d}(x_1)
+dF_{q,d}\left(\frac{1-x_1}{d}\right)
=
S_{q,d}(x_1),
\]
where we have used
\[
\frac{1-x_1}{d}=\frac{d+1-qx_1}{dq}
\]
and defined
\begin{equation}
S_{q,d}(x)
:=
F_{q,d}(x)
+dF_{q,d}\left(\frac{d+1-qx}{dq}\right),
\qquad
\frac dq\leqslant x<\frac{d+1}{q}.
\end{equation}

Equality holds in the Jensen bound when the smaller
coordinates are equal. We therefore need to prove that
$S_{q,d}$ itself is strictly positive.
We do this by showing that $S_{q,d}(d/q)>0$ and that
$S_{q,d}$ is strictly increasing.

To express its derivative, define
\begin{equation}\label{eq:Gdef}
\begin{aligned}
G_{q,d}(x)
:={}&
(qx)^{\frac{d+3}{2d}}
-\left(\frac{d+1-qx}{d}\right)^{\frac{d+3}{2d}}
+\left(\frac{d+1-qx}{d}\right)^{\frac{3d+1}{2}}\\
&-\left(1+\frac{qx-1}{d^2}\right)^{\frac{3d+1}{2}}
+\frac{d^2}{d^2+qx-1}
-\frac1{qx},
\qquad
\frac dq\leqslant x\leqslant\frac{d+1}{q}.
\end{aligned}
\end{equation}
All exponents are positive, and the denominators remain
positive on $[d/q,(d+1)/q]$. Hence $G_{q,d}$ is continuous
on this closed interval. We include the right endpoint
only in the domain of $G_{q,d}$; the domain of $S_{q,d}$
remains $[d/q,(d+1)/q)$.

\begin{lemma}[Compensation inequality]\label{lem:compensation}
For every integer $d\geqslant2$,
\[
S'_{q,d}(x)
=
\frac{3q(d+1)}{2d}\,G_{q,d}(x)
\qquad
\left(\frac dq\leqslant x<\frac{d+1}{q}\right),
\]
and
\begin{equation}\label{eq:slope}
G_{q,d}(x)>\frac1{100}
\qquad
\left(\frac dq\leqslant x\leqslant\frac{d+1}{q}\right).
\end{equation}
Consequently, $S_{q,d}$ is strictly increasing on $[d/q,(d+1)/q)$ and
\begin{equation}\label{eq:S-pos}
S_{q,d}(x)>0
\qquad(\frac dq\leqslant x<\frac{d+1}{q}).
\end{equation}
\end{lemma}

\begin{proof}
\emph{Step 1. The derivative identity.}

By the chain rule,
\[
\begin{aligned}
S_{q,d}'(x)
&=
F_{q,d}'(x)
+dF_{q,d}'\left(\frac{d+1-qx}{dq}\right)
\left(-\frac1d\right)\\
&=
F_{q,d}'(x)-F_{q,d}'\left(\frac{d+1-qx}{dq}\right).
\end{aligned}
\]
For $\frac dq\leqslant x<\frac{d+1}{q}$, the second argument belongs to
$(0,1/(dq)]$, so \eqref{eq:F1} applies at both arguments.
Moreover,
\[
\frac{d+1-q\frac{d+1-qx}{dq}}{d}
=
\frac{d^2+qx-1}{d^2}
=
1+\frac{qx-1}{d^2}.
\]
Substituting into \eqref{eq:F1} gives
\[
\begin{aligned}
\frac{2d}{3q(d+1)}S_{q,d}'(x)
={}&
\left[
(qx)^{\frac{d+3}{2d}}
+\left(\frac{d+1-qx}{d}\right)^{\frac{3d+1}{2}}
-\frac1{qx}-\frac{d}{d+1-qx}
\right]\\
&-
\left[
\left(\frac{d+1-qx}{d}\right)^{\frac{d+3}{2d}}
+\left(1+\frac{qx-1}{d^2}\right)^{\frac{3d+1}{2}}
-\frac{d}{d+1-qx}
-\frac{d^2}{d^2+qx-1}
\right]\\
={}&G_{q,d}(x).
\end{aligned}
\]
The two terms involving $d/(d+1-qx)$ cancel exactly.
It remains to prove \eqref{eq:slope}.

\medskip
\emph{Step 2. Power estimates for the endpoint bounds.}

We use the following three estimates to bound the endpoint values in Step 3.
\begin{equation}\label{eq:pow1}
d^{\frac{d+3}{2d}}\geqslant3,
\qquad
(d+1)^{\frac{d+3}{2d}}\geqslant4
\qquad(d\geqslant3),
\end{equation}
\begin{equation}\label{eq:pow2}
\left(1+\frac1d\right)^{\frac{3d+1}{2}}
<
e^{3/2}
<
\frac{449}{100}
\qquad(d\geqslant2),
\end{equation}
and
\begin{equation}\label{eq:pow3}
\left(1+\frac{d-1}{d^2}\right)^{\frac{3d+1}{2}}
<
d^{\frac{d+3}{2d}}
\qquad(d\geqslant3).
\end{equation}
We use the elementary numerical bounds
\[
\log4<\frac32,
\qquad
\log5>\frac85,
\qquad
e^{3/2}<\frac{449}{100}.
\]
In the derivative calculations in this step, $d$ is regarded
as a real variable.

For \eqref{eq:pow1}, apply $\log u\leqslant u-1$ with
$u=(d+1)/4$:
\[
\begin{aligned}
\log(d+1)
&=
\log4+\log\left(\frac{d+1}{4}\right)\\
&\leqslant
\log4+\frac{d+1}{4}-1\\
&=
\log4+\frac{d-3}{4}\\
&<
\frac32+\frac{d-3}{4}
=
\frac{d+3}{4}.
\end{aligned}
\]
Consequently,
\[
\begin{aligned}
d(d+3)-3(d+1)\log(d+1)
&>
d(d+3)-\frac{3(d+1)(d+3)}4\\
&=
\frac{(d+3)(d-3)}4
\geqslant0.
\end{aligned}
\]
It follows that
\[
\left[\frac{d+3}{2d}\log(d+1)\right]'
=
\frac{d(d+3)-3(d+1)\log(d+1)}
     {2d^2(d+1)}
>0.
\]
Thus $(d+1)^{(d+3)/(2d)}$ is increasing for $d\geqslant3$.
Its value at $d=3$ is $4$, proving the second inequality
in \eqref{eq:pow1}.

For the first inequality, note that
\[
0<\frac{d}{d+1}<1,
\qquad
0<\frac{d+3}{2d}\leqslant1.
\]
For a base in $(0,1)$, decreasing the exponent increases
the value. Hence
\[
\left(\frac{d}{d+1}\right)^{\frac{d+3}{2d}}
\geqslant\frac{d}{d+1},
\]
and therefore
\[
\begin{aligned}
d^{\frac{d+3}{2d}}
&=
(d+1)^{\frac{d+3}{2d}}
\left(\frac{d}{d+1}\right)^{\frac{d+3}{2d}}\\
&\geqslant
\frac{4d}{d+1}
=
3+\frac{d-3}{d+1}
\geqslant3.
\end{aligned}
\]

For \eqref{eq:pow2}, the identity
\[
\frac1{1+x}=1-\frac{x}{1+x}
\]
gives
\[
\begin{aligned}
\log\left(1+\frac1d\right)
&=
\int_0^{1/d}\frac1{1+x}\,dx\\
&=
\frac1d-\int_0^{1/d}\frac{x}{1+x}\,dx.
\end{aligned}
\]
On the integration interval,
\[
0\leqslant x\leqslant\frac1d
\quad\Longrightarrow\quad
1+x\leqslant\frac{d+1}{d}.
\]
Taking reciprocals and multiplying by $x\geqslant0$ gives
\[
\frac{x}{1+x}\geqslant\frac{d}{d+1}\,x.
\]
Integrating this lower bound yields
\[
\begin{aligned}
\int_0^{1/d}\frac{x}{1+x}\,dx
&\geqslant
\frac{d}{d+1}\int_0^{1/d}x\,dx\\
&=
\frac{d}{d+1}\cdot\frac1{2d^2}
=
\frac1{2d(d+1)}.
\end{aligned}
\]
Since this integral is subtracted, we obtain
\[
\log\left(1+\frac1d\right)
\leqslant
\frac1d-\frac1{2d(d+1)}.
\]
Multiplying by the positive number $(3d+1)/2$ gives
\[
\begin{aligned}
\frac{3d+1}{2}\log\left(1+\frac1d\right)
&\leqslant
\frac{3d+1}{2}
\left(\frac1d-\frac1{2d(d+1)}\right)\\
&=
\frac{(3d+1)(2d+1)}{4d(d+1)}\\
&=
\frac32-\frac{d-1}{4d(d+1)}
<
\frac32.
\end{aligned}
\]
The last inequality uses $d\geqslant2$.
Exponentiating proves \eqref{eq:pow2}.

For \eqref{eq:pow3}, direct evaluation gives the cases $d=3,4$:
\[
\left(\frac{11}{9}\right)^5<3,
\qquad
\left(\frac{19}{16}\right)^{13/2}<4^{7/8}.
\]
For real $d\geqslant5$,
\[
\left[
\log d-\frac{3d^2-2d-1}{d(d+3)}
\right]'
=
\frac{(d-1)^2(d-3)}{d^2(d+3)^2}
\geqslant0.
\]
Thus the expression in brackets is nondecreasing, and
\[
\log d-\frac{3d^2-2d-1}{d(d+3)}
\geqslant
\log5-\frac85
>0.
\]
Applying $\log(1+u)\leqslant u$ with
$u=(d-1)/d^2$, and multiplying by the positive number
$d(3d+1)/(d+3)$, gives
\[
\begin{aligned}
\frac{d(3d+1)}{d+3}
\log\left(1+\frac{d-1}{d^2}\right)
&\leqslant
\frac{d(3d+1)}{d+3}\frac{d-1}{d^2}\\
&=
\frac{3d^2-2d-1}{d(d+3)}
<
\log d.
\end{aligned}
\]
Multiplication by $(d+3)/(2d)>0$ therefore yields
\[
\frac{3d+1}{2}
\log\left(1+\frac{d-1}{d^2}\right)
<
\frac{d+3}{2d}\log d.
\]
Since the exponential function is strictly increasing,
this proves \eqref{eq:pow3}.

\medskip
\emph{Step 3. A concave lower bound when $d\geqslant3$.}

Since $0<(d+3)/(2d)\leqslant1$, the function
\[
x\longmapsto
\left(\frac{d+1-qx}{d}\right)^{\frac{d+3}{2d}}
\]
is concave on $[d/q,(d+1)/q]$.
Indeed, its second derivative on the open interval is
\[
\frac{q^2(d+3)(3-d)}{4d^4}
\left(\frac{d+1-qx}{d}\right)^{-\frac{3(d-1)}{2d}}
\leqslant0,
\]
and the function is continuous at both endpoints.

Its value and first derivative at $x=d/q$ are, respectively,
\[
d^{-\frac{d+3}{2d}},
\qquad
-\frac{q(d+3)}{2d}\,d^{-\frac{d+3}{2d}}.
\]
A concave function lies below its tangent line.
Therefore
\[
\begin{aligned}
\left(\frac{d+1-qx}{d}\right)^{\frac{d+3}{2d}}
&\leqslant
d^{-\frac{d+3}{2d}}
-\frac{q(d+3)}{2d}\,
d^{-\frac{d+3}{2d}}(x-d/q)\\
&=
d^{-\frac{d+3}{2d}}
\left(1-\frac{d+3}{2d}(qx-d)\right).
\end{aligned}
\]
This term appears with a minus sign in $G_{q,d}$, so its upper
bound gives a lower bound for $G_{q,d}$.
We may also discard the nonnegative term
$((d+1-qx)/d)^{(3d+1)/2}$.
Thus
\begin{equation}\label{eq:G-lower}
\begin{aligned}
G_{q,d}(x)\geqslant{}&
(qx)^{\frac{d+3}{2d}}
-\left(1+\frac{qx-1}{d^2}\right)^{\frac{3d+1}{2}}
+\frac{d^2}{d^2+qx-1}-\frac1{qx}\\
&-
d^{-\frac{d+3}{2d}}
\left(1-\frac{d+3}{2d}(qx-d)\right).
\end{aligned}
\end{equation}

We next show that the right-hand side is concave.
The second derivative of its first two terms is
\[
\begin{aligned}
&\left[
(qx)^{\frac{d+3}{2d}}
-\left(1+\frac{qx-1}{d^2}\right)^{\frac{3d+1}{2}}
\right]''\\
&\qquad=
\frac{q^2(d+3)(3-d)}{4d^2}
(qx)^{-\frac{3(d-1)}{2d}}
-\frac{q^2(3d+1)(3d-1)}{4d^4}
\left(1+\frac{qx-1}{d^2}\right)^{\frac{3(d-1)}2}
<0.
\end{aligned}
\]
Here the first term is nonpositive because $d\geqslant3$,
and the second is strictly negative.

For the difference of the reciprocal terms, we use
\[
d^2+qx-1-dqx
=
(d-1)(d+1-qx)
\geqslant0.
\]
Hence $d^2+qx-1\geqslant dqx>0$, which gives
\[
\frac{2q^2d^2}{(d^2+qx-1)^3}
\leqslant
\frac{2q^2d^2}{d^3(qx)^3}
=
\frac{2q^2}{d(qx)^3}.
\]

Consequently,
\[
\begin{aligned}
\left(\frac{d^2}{d^2+qx-1}-\frac1{qx}\right)''
&=
\frac{2q^2d^2}{(d^2+qx-1)^3}-\frac{2q^2}{(qx)^3}\\
&\leqslant
\frac{2q^2}{d(qx)^3}-\frac{2q^2}{(qx)^3}
=
-\frac{2q^2(d-1)}{d(qx)^3}<0.
\end{aligned}
\]
This reciprocal calculation is valid for every $d\geqslant2$
and will also be used below.
Finally, the last term in \eqref{eq:G-lower} is affine in $x$,
so its second derivative is zero.
The entire right-hand side is therefore concave.

We now bound its two endpoint values.
At the left endpoint, first note that
\[
\frac{11}{9}
-\left(1+\frac{d-1}{d^2}\right)
=
\frac{(2d-3)(d-3)}{9d^2}
\geqslant0.
\]
Both quantities being compared are positive, so taking
reciprocals gives
\[
\frac{d^2}{d^2+d-1}\geqslant\frac9{11}.
\]
Also, \eqref{eq:pow1} and $d\geqslant3$ give
\[
d^{-\frac{d+3}{2d}}\leqslant\frac13,
\qquad
\frac1d\leqslant\frac13.
\]
By \eqref{eq:pow3}, the difference of the first two powers
is strictly positive.
Hence the value of the lower bound at $x=d/q$ satisfies
\[
\begin{aligned}
&
d^{\frac{d+3}{2d}}
-\left(1+\frac{d-1}{d^2}\right)^{\frac{3d+1}{2}}
-d^{-\frac{d+3}{2d}}
+\frac{d^2}{d^2+d-1}-\frac1d\\
&\hspace{2em}>
-\frac13+\frac9{11}-\frac13
=
\frac5{33}
>
\frac1{100}.
\end{aligned}
\]

At the right endpoint,
\[
1-\frac{d+3}{2d}=\frac{d-3}{2d}\geqslant0.
\]
Combining this with $d^{-(d+3)/(2d)}\leqslant1/3$ yields
\[
-\frac{d-3}{2d}\,d^{-\frac{d+3}{2d}}
\geqslant
-\frac{d-3}{6d}.
\]
Using \eqref{eq:pow1} and \eqref{eq:pow2}, the value of the
lower bound at $x=(d+1)/q$ is therefore
\[
\begin{aligned}
&
(d+1)^{\frac{d+3}{2d}}
-\left(1+\frac1d\right)^{\frac{3d+1}{2}}
-\frac{d-3}{2d}\,d^{-\frac{d+3}{2d}}
+\frac{d-1}{d+1}\\
&\hspace{2em}>
4-\frac{449}{100}
-\frac{d-3}{6d}
+\frac{d-1}{d+1}\\
&\hspace{2em}=
\frac1{100}
+\frac{(2d-1)(d-3)}{6d(d+1)}
\geqslant
\frac1{100}.
\end{aligned}
\]
The first inequality is strict because the power in
\eqref{eq:pow2} is strictly less than $449/100$.
The final fraction is nonnegative for $d\geqslant3$.

For $x\in[d/q,(d+1)/q]$, the chord joining these two endpoint
values has weights $d+1-qx$ and $qx-d$.
These weights are nonnegative and sum to $1$.
Both endpoint values exceed $1/100$, so every value on
the chord exceeds
\[
(d+1-qx)\frac1{100}+(qx-d)\frac1{100}
=
\frac1{100}.
\]
The concave right-hand side of \eqref{eq:G-lower} lies
above this chord.
Since $G_{q,d}$ lies above that right-hand side,
\eqref{eq:slope} follows for $d\geqslant3$.

\medskip
\emph{Step 4. The case $d=2$.}

In this case, discard the nonnegative term
$((3-qx)/2)^{7/2}$ from \eqref{eq:Gdef} to obtain
\[
G_{q,2}(x)
\geqslant
(qx)^{5/4}
-\left(\frac{qx+3}{4}\right)^{7/2}
-\left(\frac{3-qx}{2}\right)^{5/4}
+\frac4{qx+3}-\frac1{qx}.
\]
We again prove that this explicit lower bound is concave
and positive at its endpoints.

For $\frac2q\leqslant x\leqslant\frac3q$,
\[
(qx)^{-3/4}\leqslant1,
\qquad
\left(\frac{qx+3}{4}\right)^{3/2}\geqslant1.
\]
The first inequality uses $qx\geqslant1$ and a negative exponent;
the second uses $(qx+3)/4\geqslant1$ and a positive exponent.

It follows that
\[
\begin{aligned}
\left[
(qx)^{5/4}
-\left(\frac{qx+3}{4}\right)^{7/2}
\right]''
&=
\frac{5q^2}{16}(qx)^{-3/4}
-\frac{35q^2}{64}
\left(\frac{qx+3}{4}\right)^{3/2}\\
&\leqslant
\frac{5q^2}{16}-\frac{35q^2}{64}
=
-\frac{15q^2}{64}<0.
\end{aligned}
\]
For $2/q<x<3/q$, we also have
\[
\left[
-\left(\frac{3-qx}{2}\right)^{5/4}
\right]''
=
-\frac{5q^2}{64}
\left(\frac{3-qx}{2}\right)^{-3/4}
<0.
\]
The calculation in Step 3, applied with $d=2$, shows that
the reciprocal difference is concave as well.
Thus the lower bound is concave on $(2/q,3/q)$.
Since all its terms are continuous on $[2/q,3/q]$,
it is concave on the closed interval.

At the two endpoints, direct evaluation gives
\[
2^{5/4}
-\left(\frac54\right)^{7/2}
-2^{-5/4}
+\frac3{10}
>
\frac1{24},
\]
and
\[
3^{5/4}
-\left(\frac32\right)^{7/2}
+\frac13
>
\frac1{24}.
\]
For $x\in[2/q,3/q]$, the chord joining these endpoint values
has weights $3-qx$ and $qx-2$. These weights are nonnegative
and sum to $1$, so the chord is strictly above $1/24$
throughout the interval. By concavity,
\[
G_{q,2}(x)>\frac1{24}>\frac1{100}
\qquad
\left(\frac2q\leqslant x\leqslant\frac3q\right).
\]
This completes the proof of \eqref{eq:slope}.

\medskip
\emph{Step 5. Positivity of the total contribution.}

By Step 1 and \eqref{eq:slope},
\[
S'_{q,d}(x)
=
\frac{3q(d+1)}{2d}G_{q,d}(x)>0
\qquad
\left(\frac dq\leqslant x<\frac{d+1}{q}\right).
\]
Hence $S_{q,d}$ is strictly increasing on $[d/q,(d+1)/q)$.

Lemma~\ref{lem:endpoint} gives $F_{q,d}(d/q)>0$.
Since $d\geqslant2$, we have $0<1/(dq)<1/q$,
so Lemma~\ref{lem:shape}(i) also gives
$F_{q,d}(1/(dq))>0$.
Hence
\[
S_{q,d}(d/q)=F_{q,d}(d/q)+dF_{q,d}(1/(dq))>0.
\]
By monotonicity,
\[
S_{q,d}(x)\geqslant S_{q,d}(d/q)>0
\qquad(\frac dq\leqslant x<\frac{d+1}{q}),
\]
which proves \eqref{eq:S-pos}.
\end{proof}

%% file: sections/05-uniform-convergence.tex
\section{From critical descent to uniform boundary convergence}
\label{sec:uniform}

Recall that
\[
p_c=p_c(d,q)
=
\max\left\{0,\,1-\frac{q}{d+1}\right\}.
\]
We now use the critical descent inequality to prove uniform
boundary convergence on $\mathbb{T}^d_n$.
All suprema below are taken over admissible boundary conditions
on $\partial\mathbb{T}^d_n$.

\subsection{Monotonicity in the edge parameter}

\begin{lemma}[Monotonicity in $p$]\label{lem:monotone}
Let $d\geqslant1$ and $q\geqslant2$ be integers.
For every $c>0$ and $y\in E_q$, the function
\[
p\longmapsto\Phi_c\bigl(T_{d,p}(y)\bigr)
\]
is nonincreasing on $[0,1]$.
\end{lemma}

\begin{proof}
Fix $c>0$ and $y\in E_q$, write $x=x(y)$, and set $a=1-p$.
Since $0<x_i<1$, all factors $1-ax_i$ are positive
for $a\in[0,1]$. Define
\[
g_i(a):=\frac{x_i}{1-ax_i},
\qquad
h_i(a):=\exp\bigl(cT_{d,1-a,i}(y)\bigr)
=
\frac{(1-ax_i)^{cd}}
     {\left(\prod_{j=1}^q(1-ax_j)\right)^{cd/q}}.
\]
Then
\[
\Phi_c\bigl(T_{d,1-a}(y)\bigr)
=
\sum_{i=1}^q h_i(a)-q.
\]
Writing $\bar g=q^{-1}\sum_{j=1}^q g_j$ and differentiating,
we obtain
\[
\frac{d}{da}T_{d,1-a,i}(y)
=
-d(g_i-\bar g),
\qquad
h_i'=-cd\,h_i(g_i-\bar g).
\]
Consequently,
\begin{equation}\label{eq:dsum-h}
\frac{d}{da}\sum_{i=1}^q h_i
=
-\frac{cd}{2q}
\sum_{i,j=1}^q(h_i-h_j)(g_i-g_j).
\end{equation}

For fixed $a$, the function $s\mapsto s/(1-as)$ is increasing
on $(0,1)$. Thus the $g_i$ have the same ordering as the $x_i$.
The denominator in the definition of $h_i$ is independent of $i$,
and $s\mapsto(1-as)^{cd}$ is nonincreasing.
Thus the $h_i$ have the reverse ordering to the $x_i$.
Consequently,
\[
(h_i-h_j)(g_i-g_j)\leqslant0
\qquad(i,j\in[q]).
\]
The derivative in \eqref{eq:dsum-h} is therefore nonnegative.
Since $a=1-p$, the required monotonicity follows.
\end{proof}

\subsection{Uniform convergence of subtree root marginals}

\begin{proposition}[Uniform boundary convergence]\label{prop:uniform}
Let $d\geqslant1$ and $q\geqslant2$ be integers, and suppose that
\[
p\in[p_c,1],
\qquad
(q,p)\neq(d+1,0).
\]
For every fixed integer $m\geqslant0$,
\begin{equation}\label{eq:unif-messages}
\lim_{n\to\infty}
\sup_{\xi}
\max_{v\in\partial\mathbb T_m^d}
\bigl\|x^{(v)}-\pi\bigr\|_\infty
=0,
\end{equation}
where $x^{(v)}$ is the root marginal of the descendant subtree at $v$, as defined in \eqref{eq:Zv}, and
\[
\pi=\left(\frac1q,\ldots,\frac1q\right).
\]
In particular, the root convergence in \eqref{eq:Delta} holds.
\end{proposition}

\begin{proof}
If $p=1$, the recursion \eqref{eq:recursion} gives
\[
x_i^{(v)}=\frac1q
\]
at every nonleaf vertex, independently of $\xi$.
Thus the quantity in \eqref{eq:unif-messages} is zero whenever $n>m$.
We henceforth assume $p<1$.

We first treat $d=1$.
In this case, each nonleaf vertex has a single child, and the recursion becomes
\[
x_i'
=
\frac{1-(1-p)x_i}{q-1+p}.
\]
Therefore
\[
x_i'-\frac1q
=
-\lambda\left(x_i-\frac1q\right),
\qquad
\lambda=\frac{1-p}{q-1+p}.
\]
Iterating from the leaf gives
\begin{equation}\label{eq:line-contraction}
\bigl\|x^{(v)}-\pi\bigr\|_\infty
\leqslant
\lambda^{\,n-m}
\qquad(|v|=m).
\end{equation}
Here $\lambda<1$ precisely when $q-2+2p>0$.
Under the assumptions of the proposition, this excludes only $(q,p)=(2,0)$, so \eqref{eq:unif-messages} follows.

It remains to consider $d\geqslant2$.
We first establish strict descent for all parameters in the proposition.
Set
\[
r=\max\{d,q-1\},
\qquad
c=c_r=\frac{3(r+1)}{2r}.
\]
Then $q\leqslant r+1$ and
\[
p_c(r,q)=p_c(d,q).
\]
Indeed, if $q\leqslant d+1$, then $r=d$.
If $q\geqslant d+2$, then $r=q-1$ and both critical parameters are zero.

By \eqref{eq:T}, the branching number appears only as a prefactor, so
\[
T_{d,p}=\frac dr\,T_{r,p}.
\]
Since $0<d/r\leqslant1$ and $\Phi_c(0)=0$, convexity, Lemma~\ref{lem:monotone}, and Theorem~\ref{thm:descent} give
\begin{equation}\label{eq:descent-p}
\begin{aligned}
\Phi_c\bigl(T_{d,p}(y)\bigr)
&\leqslant
\frac dr\,\Phi_c\bigl(T_{r,p}(y)\bigr)\\
&\leqslant
\frac dr\,\Phi_c\bigl(T_{r,p_c}(y)\bigr)\\
&<
\frac dr\,\Phi_c(y)
\leqslant
\Phi_c(y),
\qquad y\in E_q\setminus\{0\}.
\end{aligned}
\end{equation}
The map $T_{d,p}$ is continuous and fixes $0$, so it satisfies the assumptions of Lemma~\ref{lem:compact}.
Notice also that, when $q\geqslant d+2$, we have $r=q-1>d$ and
\begin{equation}\label{eq:descent-r}
\Phi_{c_r}\bigl(T_{d,p}(y)\bigr)
\leqslant
\frac dr\,\Phi_{c_r}(y)
\qquad(y\in E_q).
\end{equation}

We next obtain a uniform bound on the initial potential.
The number of layers needed for this bound depends on whether $p$ is positive.

If $0<p<1$, Lemma~\ref{lem:initial} shows that every vertex at level $n-1$ has a strictly positive root marginal and
\[
\Phi_c\bigl(y^{(v)}\bigr)
\leqslant
q\bigl(p^{-cd}-1\bigr).
\]
In this case, set
\[
\ell=1,
\qquad
a_0=q\bigl(p^{-cd}-1\bigr)>0.
\]

If $p=0$, the assumptions of the proposition imply $q\geqslant d+2$.
Every leaf boundary condition is then admissible: coloring upward from the leaves, each vertex needs to avoid at most the $d$ colors assigned to its children.

Consider the measure on the descendant subtree rooted at a nonleaf vertex $v$.

Conditionally on all other spins in this subtree, the spin at
$v$ is uniform over the colors not used by its children.
Since these $d$ children are its only neighbors in the subtree,
at least $q-d$ colors are available. Averaging over the
conditioned spins gives
\begin{equation}\label{eq:zero-temp-bound}
x_i^{(v)}\leqslant\frac1{q-d}
\qquad
(i\in[q],\ |v|\leqslant n-1).
\end{equation}

Now suppose $|v|\leqslant n-2$, so that every child of $v$
is a nonleaf vertex. Set
\[
b:=1-\frac1{q-d}\in[1/2,1).
\]
By \eqref{eq:zero-temp-bound}, every factor
$1-x_i^{(vu)}$ in the recursion belongs to $[b,1]$.
Thus $x^{(v)}$ is strictly positive, and
\[
b^d\leqslant
\frac{x_i^{(v)}}{x_j^{(v)}}
\leqslant b^{-d}
\qquad(i,j\in[q]).
\]
The argument in Lemma~\ref{lem:initial} then gives
\[
\bigl\|y^{(v)}\bigr\|_\infty
\leqslant d\log\frac1b,
\qquad
\Phi_c\bigl(y^{(v)}\bigr)
\leqslant q\bigl(b^{-cd}-1\bigr).
\]
In this case, set
\[
\ell=2,
\qquad
a_0=q\bigl(b^{-cd}-1\bigr)>0.
\]

In both cases, every vertex at level $n-\ell$ has a strictly
positive root marginal with potential at most $a_0$.
The constants $\ell$ and $a_0$ depend only on $d,q,p$.
The recursion preserves positivity, so Lemma~\ref{lem:logrec}
gives
\[
y^{(v)}
=
\frac1d\sum_{u=1}^d
T_{d,p}\bigl(y^{(vu)}\bigr)
\]
at every level closer to the root.

Fix $\varepsilon\in(0,a_0]$.
By Lemma~\ref{lem:compact}, there is an integer
$N_\varepsilon$, independent of $n$ and $\xi$, such that
\[
\Phi_c\bigl(y^{(v)}\bigr)\leqslant\varepsilon
\]
whenever $|v|=m$ and $n-m-\ell\geqslant N_\varepsilon$.
Equation~\eqref{eq:x-from-y} then gives
\[
\bigl\|x^{(v)}-\pi\bigr\|_\infty
\leqslant
\frac1q\left[
\exp\left(\frac{2(q-1)\sqrt{2\varepsilon}}c\right)-1
\right].
\]
This bound is uniform over $\xi$ and all vertices at level $m$.
Letting $\varepsilon\downarrow0$ proves
\eqref{eq:unif-messages}.
\end{proof}

\begin{remark}
The initial potential bound need only be finite.
In particular, Lemma~\ref{lem:compact} applies even when
the bound $q(p^{-cd}-1)$ is large because $p>0$ is small.

For $2\leqslant q\leqslant d$ at $p=p_c$, the compactness argument above establishes uniform convergence without computing a rate.
The precise critical rate is given by Corollary~\ref{cor:powerlaw}.
When $d\geqslant2$, $q\geqslant d+2$, and $p<1$, equation \eqref{eq:descent-r} and convexity give the geometric bound
\[
\sup_\xi
\max_{v\in\partial\mathbb T_m^d}
\Phi_{c_r}\bigl(y^{(v)}\bigr)
\leqslant
a_0\left(\frac dr\right)^{n-m-\ell},
\qquad
r=q-1,\quad n\geqslant m+\ell,
\]
with $a_0$ and $\ell$ chosen in the proof.
For $d=1$, the explicit contraction is given in \eqref{eq:line-contraction}.
\end{remark}

%% file: sections/06-gibbs-uniqueness.tex
\section{Gibbs uniqueness and the complementary parameter range}
\label{sec:gibbs}
We now complete the proof of Theorem~\ref{thm:main}.
We use the standard characterization of Gibbs uniqueness:
the $q$-state Potts model on the infinite $(d+1)$-regular tree
has a unique Gibbs measure if and only if the root marginal
on the finite rooted $d$-ary tree converges to the uniform
distribution, uniformly over admissible boundary conditions;
see \cite{GWY} and \cite[Lemma~1.3]{dBBR}.

In our notation, color symmetry allows us to write this
criterion as
\begin{equation}\label{eq:root-criterion}
\lim_{n\to\infty}
\max_{\xi}
\left|
\mu^\xi_{n,p,q}[\sigma_\ast=1]-\frac1q
\right|
=0.
\end{equation}
The same criterion characterizes uniqueness on the infinite
rooted tree $\mathbb{T}^d$.
Indeed, under \eqref{eq:root-criterion}, the factorization
in Proposition~\ref{prop:finite-subtrees} gives uniform
convergence on every fixed finite subtree.
The DLR equations then force every Gibbs measure to have
these limiting finite-dimensional marginals.

Conversely, if \eqref{eq:root-criterion} fails, there is a
sequence of admissible boundary conditions for which the
root probability of color $1$ remains a fixed positive
distance from $1/q$.
By compactness, a subsequence of the corresponding measures
converges locally. The limit satisfies the DLR equations
and differs from the free Gibbs measure at the root.

\subsection{Uniqueness and convergence on finite subtrees}

Suppose $p\in[p_c,1]$ and $(q,p)\ne(d+1,0)$.
Proposition~\ref{prop:uniform} gives \eqref{eq:root-criterion},
so the Gibbs measure is unique on both trees.
It is the free Gibbs measure $\nu$, obtained by assigning
a uniform color to the root and using the transition
probabilities $K_p$ in~(1.8) along the edges.
Its marginal on $\mathbb{T}^d_m$ is
\begin{equation}\label{eq:free-marginal}
\nu_m[\eta]
=
\frac{
\displaystyle
\prod_{\{v,w\}\in E(\mathbb{T}^d_m)}
p^{\delta_{\eta_v,\eta_w}}
}{
q(q-1+p)^{|E(\mathbb{T}^d_m)|}
},
\qquad
\eta\in[q]^{\mathbb{T}^d_m}.
\end{equation}
We next prove that the finite-volume measures converge
to these marginals uniformly over boundary conditions.

\begin{proposition}\label{prop:finite-subtrees}
Let $d\geqslant1$, $q\geqslant2$, and $p\in[p_c,1]$ with
$(q,p)\ne(d+1,0)$.
For every fixed $m\geqslant0$,
\begin{equation}\label{eq:ball-conv}
\lim_{n\to\infty}
\max_{\xi}
\max_{\eta\in[q]^{\mathbb{T}^d_m}}
\left|
\mu^\xi_{n,p,q}
\big[\sigma|_{\mathbb{T}^d_m}=\eta\big]
-
\nu_m[\eta]
\right|=0.
\end{equation}
Consequently, \eqref{eq:unifconv} holds for every event depending
only on the colors in $\mathbb{T}^d_m$.
\end{proposition}

\begin{proof}
Fix $n>m$ and an admissible boundary condition $\xi$.
For a vertex $v$ at level $m$, write $vu$ for its $u$-th child,
where $u\in[d]$.
Summing over the descendant subtrees below level $m$ gives
\begin{equation}\label{eq:ball-weight}
\omega^\xi_{n,m}(\eta)
:=
\prod_{\{v,w\}\in E(\mathbb{T}^d_m)}
p^{\delta_{\eta_v,\eta_w}}
\prod_{v\in\partial\mathbb{T}^d_m}
\prod_{u=1}^{d}
\bigl(1+(p-1)x^{(vu)}_{\eta_v}\bigr),
\end{equation}
and
\begin{equation}\label{eq:ball-normalized}
\mu^\xi_{n,p,q}
\big[\sigma|_{\mathbb{T}^d_m}=\eta\big]
=
\frac{\omega^\xi_{n,m}(\eta)}
{\displaystyle
 \sum_{\zeta\in[q]^{\mathbb{T}^d_m}}
 \omega^\xi_{n,m}(\zeta)}.
\end{equation}
Here $x^{(vu)}$ is the root marginal in the descendant subtree
rooted at $vu$, with the boundary condition induced by $\xi$.
The partition functions of these subtrees are positive and
cancel in the ratio.

Set
\[
s_p:=\frac{q-1+p}{q}>0.
\]
By Proposition~\ref{prop:uniform},
\[
1+(p-1)x_i^{(vu)}\longrightarrow s_p
\]
uniformly over $\xi,v,u$, and $i$.
Since there are $d^{m+1}$ edges leaving $\mathbb{T}^d_m$,
\[
\omega^\xi_{n,m}(\eta)
\longrightarrow
s_p^{d^{m+1}}
\prod_{\{v,w\}\in E(\mathbb{T}^d_m)}
p^{\delta_{\eta_v,\eta_w}}
\]
uniformly over $\xi$ and $\eta$.
Summing over the leaves one level at a time gives
\begin{equation}
\sum_{\eta\in[q]^{\mathbb{T}^d_m}}
\prod_{\{v,w\}\in E(\mathbb{T}^d_m)}
p^{\delta_{\eta_v,\eta_w}}
=
q(q-1+p)^{|E(\mathbb{T}^d_m)|}.
\end{equation}
Thus the denominator in~\eqref{eq:ball-normalized}
converges uniformly to a strictly positive number,
and taking the ratio proves~\eqref{eq:ball-conv}.
Since $\mathbb{T}^d_m$ has finitely many configurations,
summing over those belonging to an event $A$ proves
\eqref{eq:unifconv}.
\end{proof}

\subsection{Nonuniqueness outside the stated range}

We now recall the known nonuniqueness results for the
complementary parameter range.

At $p=0$, frozen proper colorings exist when
$2\leqslant q\leqslant d+1$; see \cite[Section~2.1]{BW}.
They give distinct Gibbs measures.
When $d=1$ and $q=2$, these are the Dirac measures
supported on the two alternating colorings.

At positive temperature, first suppose that
$3\leqslant q\leqslant d$ and
\[
0<p<1-\frac{q}{d+1}.
\]
Applying \cite[Corollary~7.5]{GSV} with $\Delta=d+1$ and $B=p$
gives a semitranslation-invariant Gibbs measure on the
infinite $(d+1)$-regular tree that is not translation invariant.
It differs from the free Gibbs measure, which is translation
invariant, so uniqueness fails.

For $q=2$ and $d\geqslant2$, nonuniqueness holds when
\[
0<p<\frac{d-1}{d+1}=1-\frac2{d+1}.
\]
This is the classical threshold for the antiferromagnetic
Ising model with zero external field;
see \cite[Section~2.2]{SST}.

These results cover every parameter outside \eqref{eq:range}.
By the characterization \eqref{eq:root-criterion},
nonuniqueness holds on both trees, and the convergence
in \eqref{eq:Delta} fails.

\begin{proof}[Proof of Theorem~\ref{thm:main}]
In the range \eqref{eq:range}, Proposition~\ref{prop:uniform}
and the characterization \eqref{eq:root-criterion} give
Gibbs uniqueness on both trees.

The unique measure is the free Gibbs measure, and
Proposition~\ref{prop:finite-subtrees} gives \eqref{eq:unifconv}.
Outside this range, the results recalled above give
nonuniqueness and failure of \eqref{eq:Delta}.
\end{proof}